\documentclass [12pt] {amsart}

\usepackage [T1] {fontenc}
\usepackage {amssymb}
\usepackage {amsmath}
\usepackage {amsthm,amsfonts}
\usepackage{tikz}
\usepackage{enumerate}
\usepackage{url}
\usepackage{mathrsfs}
\usepackage{mdwlist}
\usetikzlibrary{cd}
\usepackage{graphicx}

 \usepackage[title]{appendix}

\usetikzlibrary{arrows,decorations.pathmorphing,backgrounds,positioning,fit,petri}

\usetikzlibrary{matrix}

\usepackage{hyperref}

\DeclareMathOperator {\td} {tr.deg.}

\DeclareMathOperator {\Der} {Der}

\DeclareMathOperator {\SL} {SL}

\DeclareMathOperator {\Ann} {Ann}

\DeclareMathOperator {\rk} {rk}

\DeclareMathOperator {\C} {\mathbb{C}}

\DeclareMathOperator {\h} {\mathbb{H}}

\DeclareMathAlphabet\urwscr{U}{urwchancal}{m}{n}%
\DeclareMathAlphabet\rsfscr{U}{rsfso}{m}{n}
\DeclareMathAlphabet\euscr{U}{eus}{m}{n}
\DeclareFontEncoding{LS2}{}{}
\DeclareFontSubstitution{LS2}{stix}{m}{n}
\DeclareMathAlphabet\stixcal{LS2}{stixcal}{m} {n}

\theoremstyle {plain}

\newtheorem {theorem} {Theorem} [section]
\newtheorem {lemma} [theorem] {Lemma}
\newtheorem {proposition} [theorem] {Proposition}
\newtheorem {corollary} [theorem] {Corollary}

\newtheorem {conjecture} [theorem] {Conjecture}
\newtheorem {claim} [theorem] {Claim}

\theoremstyle {definition}
\newtheorem* {definition}{Definition}
\newtheorem* {notation} {Notation}
\newtheorem {example}[theorem]{Example}

\theoremstyle {remark}
\newtheorem {remark} [theorem] {Remark}

\numberwithin{equation}{theorem}

\calclayout

\begin {document}

\title{A Closure Operator Respecting the Modular $j$-Function}

\author{Vahagn Aslanyan}
\address{School of Mathematics, University of East Anglia, Norwich, NR4 7TJ, UK}
\email{V.Aslanyan@uea.ac.uk}
\author{Sebastian Eterovi\'c}
\address{Department of Mathematics, University of  California, Berkeley, California 94720, USA}
\email{eterovic@math.berkeley.edu}
\author{Jonathan Kirby}
\address{School of Mathematics, University of East Anglia, Norwich, NR4 7TJ, UK}
\email{Jonathan.Kirby@uea.ac.uk}
\thanks{VA and JK were supported by EPSRC grant EP/S017313/1. SE was partially supported by NSF RTG grant DMS-1646385}
\date{\today}


\maketitle

\begin{abstract}
We prove some unconditional cases of the Existential Closedness problem for the modular $j$-function. For this, we show that for any finitely generated field we can find a ``convenient'' set of generators. This is done by showing that in any field equipped with functions replicating the algebraic behaviour of the modular $j$-function and its derivatives, one can define a natural closure operator in three equivalent different ways. 
\end{abstract}

\section{Introduction}

In this paper we study some cases of the Existential Closedness problem (EC) for the modular $j$-function. In simple terms, the EC problem asks to find minimal geometric conditions that an algebraic variety $V\subset\mathbb{C}^{2n}$ should satisfy to ensure that there exists a point $\mathbf{z}\in\h^{n}$ such that 
$$(\mathbf{z},j(\mathbf{z}))\in V.$$
The results of \cite{Aslanyan-Eterovic-Kirby-Diff-EC-j} and \cite{aslanyan-kirby} inform what the conditions of EC may be (both articles consider variants of EC for $j$), and the results of \cite{Eterovic-Herrero} show general cases of varieties $V\subset\mathbb{C}^{2n}$ for which one can find such $\mathbf{z}$. One should also consider the version of the problem which includes the derivatives of $j$: what are minimal geometric conditions on a variety $V\subset\mathbb{C}^{4n}$ so that there is $\mathbf{z}\in\h^{n}$ such that 
$$(\mathbf{z},j(\mathbf{z}),j'(\mathbf{z}),j''(\mathbf{z}))\in V?$$
Furthermore, in both cases one would like to know if it is possible to find such a point so that it is generic in $V$ with respect to a given finitely generated field. These questions are natural analogues of the corresponding statement for the complex exponential function (called the strong exponential closedness conjecture), first considered by Zilber in his work on pseudo-exponentiation (\cite{Zilb-pseudoexp}), and further studied in \cite{Bays-Kirby-exp,daquino-fornasiero-terzo,mantova,marker}. 

In our paper, we study the EC problem in the case when $V$ is a plane curve. This question was already considered in \cite[Theorem 1.2]{Eterovic-Herrero} where it was proven that, if one assumes a modular version of Schanuel's conjecture (which is considered out of reach), then for any irreducible plane curve $V$ which is not a horizontal or a vertical line and any finitely generated field $F$, $V$ has a point of the form $(z,j(z))$ which is generic over $F$. Our main result shows cases in which one can remove the dependence on the modular Schanuel conjecture to obtain an unconditional result. For this we use the Ax-Schanuel theorem for $j$ established in \cite{Pila-Tsim-Ax-j}. The Ax-Schanuel theorem for $j$ shows (among other things) that the modular version of Schanuel's conjecture holds when we consider elements which are ``transcendental with respect to the $j$-function''. More precisely, there exists a countable algebraically closed subfield $C\subset\mathbb{C}$ (which is not finitely generated) consisting of all the complex numbers which are considered to be algebraic with respect to $j$. We now state our main result. 
\begin{theorem} \label{thm:unconditional}
Let $V\subset\C^{2}$ be an algebraic curve which is neither a horizontal nor a vertical line. If $V$ is not definable over $C$, then for any finitely generated field $F$, there exists $z\in\h$ such that $(z,j(z))\in V$ and $\mathrm{tr.deg.}_{F}F(z,j(z)) \geq 1$, i.e.\ $(z,j(z))$ is generic in $V$ over $F$. 
\end{theorem}
An algebraic description of the elements of $C$ is obtained in \textsection\ref{subsec:khovanskii}. We will show in Remark \ref{rem:notinC} that, assuming the generalised period conjecture of Grothendieck-Andr\'e (which implies the modular version of Schanuel's conjecture, see \textsection\ref{subsec:msc}), one can show that $\pi\notin C$.  In general, finding explicit complex numbers which are not in $C$ is likely to be of similar difficulty to proving a significant part of the generalised period conjecture. However, this does depend on the definition of being ``explicit''. The methods of \cite{jones-servi} give computable complex numbers which are not exponentially algebraic, via a diagonalisation argument, and it is likely that a similar argument would give computable complex numbers which are not in $C$.

The ideas used in the proof of Theorem \ref{thm:unconditional} can also be used to get higher dimensional cases of EC. In \textsection \ref{subsec:higherdimec} we give the following example: let $j_{1}(z):=j(z)$ and for $k>1$ set $j_{k+1}(z):=j(j_{k}(z))$, then for any positive integer $n$, for any $a\notin C$ and for any finitely generated field $K$ there exists $z\in\h$ such that $j_{n}(z)$ exists, $z = j_{n}(z)+a$ and
\begin{equation*}
    \mathrm{tr.deg.}_{K}K(z,j_{1}(z),j_{2}(z),\ldots,j_{n}(z))\geq n.
\end{equation*}

The proof of Theorem \ref{thm:unconditional} relies on some calculations which we believe are interesting in their own right, regarding how to find a ``convenient'' set of generators for a given finitely generated field (see Theorem \ref{thm:mainJ}). This set of generators will consist of elements which are not in $C$. 

Furthermore, our methods can be generalised in a straightforward way to the complex exponential function to get an unconditional statement analogous to Theorem \ref{thm:unconditional}: if $E\subset\mathbb{C}$ denotes the countable algebraically closed subfield of exponentially algebraic numbers, and $V\subset\C^{2}$ is a plane curve not definable over $E$ which is neither a horizontal nor a vertical line, then for any finitely generated field $F$ over which $V$ is defined, there exists $z\in\C$ such that $(z,\exp(z))\in V$ and $\mathrm{tr.deg.}_{F}F(z,\exp(z)) = 1$, see Theorem \ref{thm:unconditionalfore} (this gives unconditional cases of the main result of \cite{mantova}, which assumes Schanuel's conjecture). We remark that similar ideas to the ones used in the proof of Theorem \ref{thm:unconditionalfore} can be found in \cite[\textsection 11]{Bays-Kirby-exp}.

Note that in Theorem \ref{thm:unconditional} we do not deal with varieties defined over $C$ as this allows us to use functional transcendence results. For varieties defined over $C$, which is an arithmetic set in nature, we need arithmetic transcendence statements (such as the modular Schanuel conjecture). For example, if $V\subset\mathbb{C}^{2}$ is defined by the equation $X=Y$ and we could apply the conclusion of Theorem \ref{thm:unconditional} to $V$, then we can show that $\td_{\mathbb{Q}} \mathbb{Q}(\{z\in \h: z=j(z) \})$ is infinite, which is an unknown transcendence statement.

\subsection{The general setting}
To prove our main result, we will work in the general setting of $j$-fields, which gives a convenient framework in which to define some pregeometries that help us study the algebraic properties of the complex $j$-function. A $j$-field consists of a field $K$ of characteristic zero endowed with some functions whose behaviour mimics that of the complex $j$-function and its derivatives on $\mathbb{C}$. We also distinguish a subset $D\subset K$ to be the domain of these functions. On these fields, one can define a natural notion of algebraicity and transcendence over the $j$-function using any of three notions of closure: a closure with respect to field derivations on $K$ respecting the $j$-function (denoted $j\mathrm{cl}$), a closure with respect to Khovanskii systems of polynomial equations involving the $j$-function (denoted $\mathrm{kcl}$), and a closure with respect to a predimension function $\delta_{j}$ arising from a modular Schanuel property (denoted $\mathrm{cl}_{\delta_{j}}$). We will prove that all these notions define in fact the same closure operator, and so for Theorem \ref{thm:unconditional} we will use $C = j\mathrm{cl}(\emptyset)$. 
\begin{theorem}
\label{thm:main}
Let $(K,D)$ be a $j$-field. For every subset $A\subseteq K$ we have that 
$$j\mathrm{cl}(A) = \mathrm{kcl}(A) = \mathrm{cl}_{\delta_{j}}(A).$$
Moreover, this operator is a pregeometry.
\end{theorem}
This theorem will be obtained as a combination of Theorems \ref{thm:predim} and \ref{thm:extender3}. The definitions of these operators are given later on. With this result we will be able to prove the existence of the convenient set of generators we mentioned earlier required in our proof of Theorem \ref{thm:unconditional}. We also remark that it was previously shown in \cite{Eterovic-Schan-for-j} that $j\mathrm{cl}$ is a pregeometry. 

Previously, an analogous equivalence of closure operators was obtained in \cite{Kirby-Schanuel} for exponential fields. A first study of $j$-fields was given in \cite{Eterovic-Schan-for-j}, however, it was left as an open question whether Theorem \ref{thm:main} could be achieved. In technical terms, the step that was missing was: show that if $(K_{0},D_{0})\hookrightarrow (K,D)$ is a self-sufficient embedding of $j$-fields, then any $j$-derivation on $K_{0}$ extends to a $j$-derivation on $K$. In this paper we show how one can adapt an argument from \cite{Aslanyan-Eterovic-Kirby-Diff-EC-j} to solve this step (see Proposition \ref{prop:extender2}).

\subsection*{Structure of the paper}

\begin{enumerate}
    \item[\textsection \ref{sec:basicj}:] We review some basic definitions and notation we will use regarding the $j$ function. 
    \item[\textsection \ref{sec:reviewofjfields}:] We recall the definition of $j$-fields and some basic properties. We also recall some properties of $j$-derivations and define the pregeometry $j\mathrm{cl}$ and the dimension $\dim^{j}$. 
    \item[\textsection \ref{sec:fingenjfields}:] We define the predimension $\delta_{j}$ and study self-sufficient embeddings of $j$-fields. At the end of this section we prove Proposition \ref{prop:extender2} (extension of $j$-derivations), we define $\mathrm{cl}_{\delta_{j}}$, and prove one half of Theorem \ref{thm:main}: Theorem \ref{thm:predim}. 
    \item[\textsection \ref{sec:genJfinfields}:] We prove Theorems \ref{thm:unconditional} and \ref{thm:mainJ}. We also show how to get similar results for the complex exponential function and prove Theorem \ref{thm:unconditionalfore}.
    \item[\textsection \ref{sec:jpolys}:] We start by introducing the notion of $j$-polynomials and their corresponding Khovanskii systems, we define the closure operator $\mathrm{kcl}$, and prove the other half of Theorem \ref{thm:main}: Theorem \ref{thm:extender3}. We give some examples where we combine Khovanskii systems with the results of \textsection\ref{sec:genJfinfields} to get higher dimensional cases of EC. We then give an overview of the modular version of Schanuel's conjecture. We finish by showing how this conjecture can be used along with Theorem \ref{thm:extender3} to obtain (conditionally) convenient sets of generators like in Theorem \ref{thm:mainJ}, but this time using elements of $C$.
\end{enumerate}
Throughout, we will freely use terminology associated with pregeometries, and we refer the reader to \cite[\textsection 8.1]{marker} or \cite[Appendix C]{tent-ziegler}.

\section{The \texorpdfstring{$j$}{j}-function}
\label{sec:basicj}
We denote by $\mathbb{H}^{+}$ the complex upper-half plane $\{z\in \mathbb{C}:\mathrm{Im}(z)>0\}$. The group $\mathrm{GL}_2^+(\mathbb{R})$ of $2\times 2$ matrices with entries in $\mathbb{R}$ and positive determinant acts on $\mathbb{H}^{+}$ via the formula 
$$gz:=\frac{az+b}{cz+d} \ \text{ for } \ g=\left(\begin{array}{cc}a & b\\ c & d\end{array}\right) \in \mathrm{GL}^+_{2}(\mathbb{R}).$$ 
Consider the group $\mathrm{GL}_2^+(\mathbb{Q})$ consisting of the elements of $\mathrm{GL}_2^+(\mathbb{R})$ with entries in $\mathbb{Q}$. A subgroup of $\mathrm{GL}_2^+(\mathbb{Q})$ is the modular group is $\mathrm{SL}_2(\mathbb{Z})$. The modular $j$-function is defined as the unique holomorphic function $j:\mathbb{H}^{+}\to \mathbb{C}$ that satisfies 
$$j(gz)=j(z) \text{ for every }g \in \mathrm{SL}_2(\mathbb{Z}) \text{ and every }z\in\mathbb{H}^{+},$$ 
and has a Fourier expansion of the form 
\begin{equation}\label{eq:j-fourier-expansion}
j(z)=q^{-1}+744 +\sum_{k=1}^{\infty}a_kq^k \text{ with }q:=\exp(2\pi i z)  \text{ and }a_k\in \mathbb{C}.
\end{equation} 
It induces an analytic isomorphism of Riemann surfaces $\mathrm{SL}_2(\mathbb{Z}) \backslash \mathbb{H}^{+}\simeq \mathbb{C}$. The quotient space $Y=\mathrm{SL}_2(\mathbb{Z})\backslash \mathbb{H}^{+}$ is known to be a moduli space for complex tori, or equivalently, elliptic curves over $\mathbb{C}$. If $\mathrm{SL}_2(\mathbb{Z}) z$ is a point in $Y$ and $E_z$ denotes an elliptic curve in the corresponding isomorphism class, then $j(z)$ is simply the $j$-invariant of the curve $E_z$. 

\subsection{Notation}
Let $G := \mathrm{GL}_{2}(\mathbb{Q})$ and let $\mathbb{H}:= \mathbb{H}^{+}\cup\mathbb{H}^{-}$ be the union of the upper and lower half-planes. As we detailed above, the $j$-function is understood as a modular function defined on $\h^{+}$, but we will extend it to be defined on all of $\h$ using Schwarz reflection, so that $j:\mathbb{H}\rightarrow\mathbb{C}$. This means that given $z\in\mathbb{H}^{-}$, we define $j(z):=\overline{j(\overline{z})}$, where $\overline{w}$ is the complex conjugate of $w$. We have extended the $j$-function in this way because the condition for an element $g\in \mathrm{GL}_{2}(\mathbb{Q})$ to preserve $\h$ setwise can be checked only using field operations: namely $\det(g)\neq 0$; whereas preserving $\h^{+}$ would require us to introduce an order relation in the structure of $j$-fields.

If $F$ is a field, we will use $\overline{F}$ to denote an algebraic closure of $F$. In most cases, we will only consider the case of $F\subset\C$, and so $\overline{F}$ means the algebraic closure inside of $\C$. 

Tuples of elements will be denoted with boldface letters, that is, if $x_{1},\ldots,x_{m}$ are elements of a set $X$, then we write $\mathbf{x}:=(x_{1},\ldots,x_{m})$ for the ordered tuple. We will also sometimes use $\mathbf{x}$ to denote the set $\left\{x_{1},\ldots,x_{m}\right\}$, which should not lead to confusion. If $f$ denotes a function defined on $X$, then we write $f(\mathbf{x})$ to mean $(f(x_{1}),\ldots,f(x_{m}))$. Furthermore, we use $J$ to denote the triple of functions $(j,j',j'')$, so that if $z_{1},\ldots,z_{n}$ are elements of $\h$, then:
$$J(\mathbf{z}) := (j(z_{1}),\ldots,j(z_{n}),j'(z_{1}),\ldots,j'(z_{n}),j''(z_{1}),\ldots,j''(z_{n})).$$

\subsection{Special points}

A point $z$ in $\mathbb{H}$ is said to be \emph{special} if there is a non-scalar matrix $g\in G$ (in other words, $g$ is a non-central element of $G$) such that $z$ is a fixed point of $g$. This is equivalent to saying that $z$ satisfies a non trivial quadratic equation with integer coefficients. By a theorem of Schneider (\cite{schneider}) we know that $\mathrm{tr.deg.}_{\mathbb{Q}}(z,j(z)) = 0$ if and only if $z$ is special. The special points of $\mathbb{H}$ are exactly those points for which the corresponding elliptic curve (more precisely, any representative in the corresponding isomorphism class of elliptic curves) has complex multiplication.  For this reason, special points are also known as \emph{CM points} in the literature. $\Sigma\subset\mathbb{H}$ will denote the set of all special points.

\subsection{Modular polynomials}\label{subsec:mod-poly} Let $\left\{\Phi_{N}(X,Y)\right\}_{N=1}^{\infty}\subseteq\mathbb{Z}[X,Y]$ denote the family of \emph{modular polynomials} associated to $j$ (see \cite[Chapter 5, Section 2]{lang} for the definition and main properties of this family). We recall that $\Phi_{N}(X,Y)$ is irreducible in $\mathbb{C}[X,Y]$, $\Phi_{1}(X,Y) = X-Y$, and for $N\geq 2$, $\Phi_{N}(X,Y)$ is symmetric of total degree $\geq 2N$. Also, the action of $G$ on $\mathbb{H}$ can be traced by using modular polynomials in the following way: for every $g$ in $G$ we define $\mathrm{red}(g)$ as the unique matrix of the form $rg$ with $r\in \mathbb{Q},r>0$ such that the entries of $rg$ are all integers and relatively prime. Then, for every $z_1,z_2$ in $\mathbb{H}$ the following statements are equivalent: 
\begin{itemize}
    \item[(M1):] $\Phi_{N}(j(z_1),j(z_2)) = 0$,
    \item[(M2):] $gz_1=z_2$ for some $g$ in $G$ with $\det\left(\mathrm{red}(g)\right) = N$.
\end{itemize}

\subsection{Ax-Schanuel for the \texorpdfstring{$j$}{j}-function}

The $j$-function satisfies an order $3$ algebraic differential equation over $\mathbb{Q}$, and none of lower order (i.e.\ its differential rank over $\mathbb{C}$ is $3$). Namely, $\Psi(j,j',j'',j''')=0$ where 
$$\Psi(y_0,y_1,y_2,y_3)=\frac{y_3}{y_1}-\frac{3}{2}\left( \frac{y_2}{y_1} \right)^2 + \frac{y_0^2-1968y_0+2654208}{2y_0^2(y_0-1728)^2}\cdot y_1^2.$$
It is well known that any function of the form $j(gz)$ with $g \in \SL_2(\mathbb{C})$ satisfies the differential equation $\Psi(y,y',y'',y''')=0$ and all solutions (not necessarily defined on $\mathbb{H}^{+}$) are of that form (see, for example, \cite[Lemma 4.2]{Freitag-Scanlon} or \cite[Lemma 4.1]{Aslanyan-adequate-predim}). 

Observe that if $z\in\mathbb{H}$ is such that $j'(z)\neq 0$ and $j(z)\neq 0, 1728$, then we can write $j'''(z)  = \eta(j(z),j'(z),j''(z))$, where:
$$\eta(y_0,y_1,y_2):= \frac{3}{2}\cdot \frac{y_2^{2}}{y_1} - \frac{y_0^2-1968y_0+2654208}{2y_0^2(y_0-1728)^2}\cdot y_1^3.$$

\begin{definition}
Let $K$ be a field of charactersitic zero. We say that two elements $x,y\in K$ are \emph{modularly independent}\footnote{Note that this notion of independence defines a pregeometry of trivial type on $K$.} if for every modular polynomial $\Phi_{N}(X,Y)$ we have that $\Phi_{N}(x,y)\neq 0$.
\end{definition}

We now recall the statement of the Ax-Schanuel theorem for the $j$ function. The version we will use is stated in terms of differential fields. Regarding notation, we should clarify that the symbols $j_{i},j_{i}',j_{i}'',j_{i}'''$ are supposed to represent abstract elements of the differential field, and it is not a priori the case that $j_{i}'$ is the derivative of $j_{i}$, etc. For abstract differential fields we will not use the symbol $'$ to denote derivatives. 

\begin{theorem}[Ax-Schanuel for $j$, {{\cite[Theorem 1.3]{Pila-Tsim-Ax-j}}}]\label{j-chapter-Ax-for-j}
Let $(K;+,\cdot,\partial_{1},\ldots,\partial_{m},0,1)$ be a differential field, where $\partial_{1},\ldots,\partial_{m}$ are commuting derivations, let $C = \ker_{k=1}^{m}\partial_{k}$, and let $z_i, j_i, j_i', j_i'', j_i''' \in K^{\times},~ i=1,\ldots,n,$ be such that for all $i\in\left\{1,\ldots,n\right\}$ and all $k\in\left\{1,\ldots,m\right\}$:
\begin{equation*}
\Psi\left(j_i,j_i', j_i'', j_i'''\right)=0 \wedge j_{i}\notin C\wedge \partial_{k} j_i=j_i'\partial_{k} z_i \wedge \partial_{k} j_i'=j_i''\partial_{k} z_i\wedge \partial_{k} j_i''=j_i'''\partial_{k} z_i.
\end{equation*}
If the $j_i$'s are pairwise modularly independent then 
\begin{equation}\label{j-chapter-Ax-ineq}
\td_CC\left(\mathbf{z},\mathbf{j},\mathbf{j}',\mathbf{j}''\right) \geq 3n+\mathrm{rank}(\partial_{k}z_{i})_{i,k}.
\end{equation}
\end{theorem}

\section{Review of \texorpdfstring{$j$}{j}-Fields}
\label{sec:reviewofjfields}
The definition of $j$-fields and some of their basic properties were laid out in \cite{Eterovic-Schan-for-j}. For the convenience of the reader, we recall them here. 

Given a field $K$ of characteristic zero, there is a natural action of $\mathrm{GL}_{2}(\mathbb{Q})$ on $\mathbb{P}^{1}(K) = K\cup\left\{\infty\right\}$, given by:
\begin{equation*}
    gx = \frac{ax+b}{cx+d},
\end{equation*} where $g\in\mathrm{GL}_{2}(\mathbb{Q})$ is represented by $g=\left(\begin{matrix}
a & b\\
c & d
\end{matrix}\right)$. Whenever we say that $\mathrm{GL}_{2}(\mathbb{Q})$ acts on $K$, it will be in this manner. Throughout, let $G = \mathrm{GL}_{2}(\mathbb{Q})$. 

\begin{definition}
Let $K$ be a field of characteristic 0. Given a subset $A$ of $K$ we define the \emph{$G$-closure of $A$}, denoted $G\mathrm{cl}(A)$, as the set of $x\in K$ such that there exist $a\in A$ and $g\in G$ satisfying $x = ga$ (which, save for the exclusion of the point at infinity, is the union of the $G$-orbits of points in $A$).\footnote{In \cite{Eterovic-Schan-for-j} this was called ``geodesic closure'', and was denoted as $\mathrm{gcl}$. We have opted to maintain the notation used in \cite{Eterovic-Herrero}.}

It is straightforward to check that $G\mathrm{cl}$ is a pregeometry. Given $A,B\subseteq K$ let $\dim_{G}(A|B)$ be the dimension defined by the pregeometry $G\mathrm{cl}$, that is, $\dim_{G}(A|B)$ is the number of distinct orbits of elements in $A$ that do not contain elements of $B$. If $B=\emptyset$ then we write simply $\dim_{G}(A)$. 
\end{definition}

\begin{definition}
A \emph{$j$-field}\footnote{This definition of $j$-fields is slightly simpler than the one presented in \cite{Eterovic-Schan-for-j}, but it defines the same structures.} is a structure $\left<\mathbb{K}; D,j,j',j'',j'''\right>$, where:
\begin{itemize}
    \item $\mathbb{K} = \left<K; +, \cdot, 0, 1\right>$ is a field of characteristic zero,
    \item $D$ is a subset of $K$,
    \item $j,j',j'',j''':D\rightarrow K$ are functions,
\end{itemize} that satisfies:
\begin{enumerate}[(a)]
\item\label{ax:D} $D$ is closed under the action of $G$.  
\item\label{ax:diffj} For every $z\in D$, 
\begin{equation*}
    (j(z)\neq 0\wedge j(z)\neq1728 \wedge j'(z)\neq 0)\implies\Psi(j(z),j'(z),j''(z),j'''(z))=0.
\end{equation*}
\item\label{ax:modpoly1} The axiom scheme: for every $z_{1},z_{2}\in D$, if $z_{1} = gz_{2}$, then $\Phi_{N}(j(z_{1}),j(z_{2}))=0$, where $N=\det(\mathrm{red}(g))$. We also include here the expressions that can be obtained by differentiating modular relations up to $3$-times. This means the following: choosing $g$ and $N$ as before, we have that for every $z\in \h$, $\Phi_{N}(j(z),j(gz))=0$. If we interpret this expression in $\mathbb{C}$ and derive it with respect to $z$, then we get: 
\begin{equation}
\label{eq:j'}
    \Phi_{N1}(j(z),j(gz))j'(z) + \Phi_{N2}(j(z),j(gz))j'(gz)\frac{ad-bc}{(cz+d)^{2}} = 0,
\end{equation}
where $\Phi_{N1}$ and $\Phi_{N2}$ are the derivatives of $\Phi_{N}(X,Y)$ with respect to the variables $X$ and $Y$ respectively. So, for each $g\in G$, we also include the axiom that says: if $z_{1}=gz_{2}$, then:
\begin{equation*}
    \Phi_{N1}(j(z_{1}),j(z_{2}))j'(z_{1}) + \Phi_{N2}(j(z_{1}),j(z_{2}))j'(z_{2})\frac{ad-bc}{(c\alpha(z_{1})+d)^{2}} =0. 
\end{equation*}
Deriving equation (\ref{eq:j'}) again with respect to $z$, we get another equation, this time involving $j''(z)$ and $j''(gz)$, that we restate as a first-order axiom as we just did with equation (\ref{eq:j'}). Deriving (\ref{eq:j'}) twice with respect to $z$, we get an equation involving $j'''(z)$ and $j'''(gz)$, which we likewise restate as an axiom. 
\item\label{ax:modpoly} The axiom scheme (one statement for each $N\in\mathbb{N}$): for all $z_{1},z_{2}\in D$, 
\begin{equation*}
    \Phi_{N}(j(z_{1}), j(z_{2})) = 0\implies \bigvee_{g\in G, \det(\mathrm{red}(g))=N} (gz_{2} = z_{1}).
\end{equation*}
This axiom is a converse to part of axiom (\ref{ax:modpoly1}).
\item\label{ax:kernel} The statement: let $A\subseteq G$ be the set of non-scalar matrices. Then for every $z\in D$, 
\begin{equation*}
    (j(z)=0\vee j(z)=1728\vee j'(z)=0)\implies \bigvee_{g\in A}(gz=z). 
\end{equation*}
Note that the values chosen are the same as those used in axiom (\ref{ax:diffj}), and correspond to the points $z$ where the expression $\Psi(j(z),j'(z),j''(z),j'''(z))$ is not defined.
\end{enumerate}
If the function $j:D\rightarrow K$ happens to be surjective (like in the case of $j:\h\rightarrow\C$), then we say that the $j$-field is \emph{full}. 
\end{definition}

Note that as $D$ is closed under the action of $G$, we get that $D\cap\mathbb{Q}=\emptyset$ because for every $x\in\mathbb{Q}$ there is $g\in G$ such that $gx=\infty$. 

Axioms (\ref{ax:modpoly1}) and (\ref{ax:modpoly}) ensure that the equivalence between (M1) and (M2) (see \textsection\ref{subsec:mod-poly}) holds on every $j$-field. In particular, as the first modular polynomial is $\Phi_{1}(X,Y) = X-Y$, then axiom (\ref{ax:modpoly1}) implies that in every $j$-field, $j$ is invariant under $\mathrm{SL}_{2}(\mathbb{Z})$. 

Axioms (a), (b), (c) are first-order expressible in the language $\mathcal{L}_{J}:=\left\{0,1, +,\cdot,D,j,j',j'',j'''\right\}$. Axioms (d) and (e) on the other hand, are not, as they require a countably infinite disjunction of statements (they are $\mathcal{L}_{\omega_{1},\omega}$-expressible though). The first-order axioms (a), (b) and (c) do not give a complete axiomatisation of the first-order theory of the $\mathcal{L}_{J}$-structure $\left<\C, \h, j, j', j'', j'''\right>$. The presence of $j'''$ is superfluous (due to the differential equation $\Psi(j,j'j,'',j''')=0$ being satisfied); we have chosen to include it just for expository reasons. We have not stated what the values of $j''(z)$ and $j'''(z)$ are when $j(z)=0$, $j(z)=1728$ or $j'(z)=0$. This is because the actual values that these functions take for such $z$ will not be relevant.

\begin{notation}
We normally denote $j$-fields simply as $(K,D)$. 
\end{notation}

\begin{definition}
Let $(K,D)$ be a $j$-field. A point $z\in D$ will be called \emph{special} if there exists a positive integer $N>1$ such that $\Phi_{N}(j(z),j(z))=0$. Let $\Sigma_{D}$ denote the set of special points in $D$.
\end{definition}

\begin{lemma}
\label{lem:sppts}
Let $(K,D)$ be a $j$-field. Then $z\in\Sigma_{D}$ if and only if there exists a non-scalar $g\in G$ such that $gz=z$. 
\end{lemma}
\begin{proof}
We use the equivalence between (M1) and (M2) (ensured by axioms (\ref{ax:modpoly1}) and (\ref{ax:modpoly}) in the definition of $j$-fields). If $z$ is special, then there exists an integer $N>1$ such that $\Phi_{N}(j(z),j(z))=0$, which means that there exists $g\in G$ with $N = \det(\mathrm{red}(g)$ such that $gz=z$. By the definition of $\mathrm{red}(g)$, we deduce that $g$ cannot be scalar. 

Now suppose there exists a non-scalar $g\in G$ such that $gz=z$. Let $N = \det(\mathrm{red}(g))$, if $N>1$, then we are done. If instead $N=1$, then it is a simple exercise to find $h\in G$ such that $hz=z$ and $\det(\mathrm{red}(h))>1$. 
\end{proof}

\begin{definition}
A \emph{morphism of $j$-fields} $\sigma:(K_{1},D_{1})\rightarrow(K_{2},D_{2})$ is a field morphism $\sigma:K_{1}\rightarrow K_{2}$ such that $\sigma(D_{1})\subset D_{2}$ and for every $f\in\left\{j,j',j'',j'''\right\}$ and every $z\in D_{1}$ we have $\sigma(f(z)) = f(\sigma(z))$. Note that field morphisms respect the action of the group $G$ on $K$, so we also have $\sigma(gz)=g\sigma(z)$ for every $g\in G$ and $z\in D_{1}$. 

A \emph{$j$-subfield} of $(K,D)$ is a $j$-field $(K_{0},D_{0})$ such that $K_{0}\subset K$ and the inclusion map $\mathrm{id}:(K_{0},D_{0})\rightarrow(K,D)$ is a morphism of $j$-fields. 
\end{definition}

\begin{definition}
If $(K_{i},D_{i})_{i\in I}$ is a family of $j$-subfields of $(K,D)$, then we define their intersection as
$$\bigcap_{i\in I}(K_{i},D_{i}) := \left(\bigcap_{i\in I} K_{i}, \bigcap_{i\in I} D_{i}\right),$$ 
so that it is again a $j$-subfield of $(K,D)$. 

If $(A,D_{A})$ and $(B,D_{B})$ are $j$-subfield of $(K,D)$, we define their union to be:
$$(A,D_{A})\cup (B,D_{B}):= (AB, D_{A}\cup D_{B}),$$
where $AB$ denotes the subfield of $K$ generated by $A\cup B$. In this way, $(A,D_{A})\cup (B,D_{B})$ is again a $j$-subfield of $(K,D)$. 
\end{definition}

\begin{definition}
Given a subset $X\subset K$ we can define the \emph{$j$-subfield generated by $X$}, which we denote as $\left<X\right>_{j}$, as the intersection of all $j$-subfields $(K',D')$ of $(K,D)$ such that $X\subseteq K'$ and $X\cap D\subseteq D'$. 

More explicitly, $\left<X\right>_{j}$ is the $j$-field $(A,D_{A})$, where $D_{A} = G\mathrm{cl}(X\cap D)$ and $A$ is the subfield of $K$ generated by $X\cup J(D_{A})$. 

If $(K',D')$ is a $j$-subfield of $(K,D)$, then we define the $j$-subfield generated by $X$ \emph{over} $(K',D')$, to be $\left<X\right>_{j}\cup(K',D')$. 

A $j$-subfield of $(K,D)$ will be called \emph{finitely generated} if it can be generated by a finite set. Given a $j$-subfield $(K',D')$, we say that the $j$-subfield $(K_{1},D_{1})$ of $(K,D)$ is \emph{finitely generated over $(K',D')$} if there exists a finite set $X\subset K$ such that $(K_{1},D_{1}) = \left<X|(K',D')\right>_{j}$.
\end{definition}

Any field $F$ of characteristic 0 can be made into a $j$-field trivially by setting $D=\emptyset$ and then having the maps $j, j', j'', j'''$ be the empty map. To prevent trivial situations like this, we have chosen to define  $\left<X\right>_{j}$ in a way that maximises the domain. 

\begin{definition}
We will say that a $j$-field $(K,D)$ is \emph{graph-generated by $J$} if $K$ is generated as a field by the set $D\cup J(D)$. 
\end{definition}

\subsection{\texorpdfstring{$j$}{j}-derivations}
\label{subsec:jderivations}
\begin{definition}
Let $K$ be a field of characteristic 0 and let $M$ be a $K$-vector space. A map $\partial:K\rightarrow M$ is a called a \emph{derivation} if it satisfies the following two conditions:
\begin{enumerate}
\item $\partial(a+b) = \partial(a)+\partial(b)$ for every $a,b\in K$.
\item $\partial(ab) = a\partial(b)+b\partial(a)$ for every $a,b\in K$.
\end{enumerate}
Given a subset $X\subseteq K$, let $\mathrm{Der}(K/X; M)$ denote the set of derivations $\partial:K\rightarrow M$ such that $X\subseteq\ker\partial$. When $M=K$ we simplify the notation to be $\mathrm{Der}(K/X):= \mathrm{Der}(K/X;K)$. Define $\Omega(K/X)$ as the $K$-vector space generated by formal symbols of the form $dr$, where $r\in K$, quotiented by the relations given by the axioms of derivations plus that for every $x\in X$, $dx=0$. Denote by $d:K\rightarrow\Omega(K/X)$ the map $r\mapsto dr$. The map $d$ is called the \emph{universal derivation on $X$}.
\end{definition}

\begin{remark}
\label{rem:cdder}
It is well-known that for any $\partial\in \mathrm{Der}(K/X)$ there exists a $K$-linear map $\partial^{\ast}$ giving a commutative diagram:
   \begin{center}
    \begin{tikzcd}[ampersand replacement=\&]
    K \arrow[d,"\partial"]\arrow[r,"d"] \& \Omega(K/X) \arrow[dl,"\partial^{\ast}"]\\
    K \& 
    \end{tikzcd} 
    \end{center}
   which allows us to identify $\mathrm{Der}(K/X)$ with the dual of $\Omega(K/X)$. 
   
    If $\partial$ is a derivation on $K$ and $u$ is algebraic over $K$, then there is a unique way of extending $\partial$ to $K(u)$ (recall we assume $\mathrm{char}(K)=0$). On the other hand, if $t$ is transcendental over $K$, then we can extend $\partial$ to $K(t)$ by choosing any value in $K(t)$ for $\partial(t)$. For this reason, if $F$ is a subfield of $K$, then 
   $$\dim\Omega(K/F) = \dim\mathrm{Der}(K/F) = \mathrm{tr.deg.}_{F}K.$$
\end{remark}

\begin{definition}
Let $(K,D)$ be a $j$-field and let $M$ be a $K$-vector space. A derivation $\partial:K\rightarrow M$ is called a \emph{$j$-derivation} if it satisfies: $\partial(j(z)) = j'(z)\partial(z)$,\quad $\partial(j'(z)) = j''(z)\partial(z)$,\quad $\partial(j''(z)) = j'''(z)\partial(z)$, for every $z\in D$.

Let $X$ be a subset of $K$. We define $j\mathrm{Der}(K/X; M)$ as the set of $j$-derivations $\partial:K\rightarrow M$ satisfying $X\subseteq\ker\partial$. For convenience we write $j\mathrm{Der}(K/X) := j\mathrm{Der}(K/X;K)$. Note that all these spaces are $K$-vector spaces. 

Let $\Xi(K/X)$ be the vector space obtained from $\Omega(K/X)$ by taking the quotient with the subspace generated by the axioms for $j$-derivations. This induces a map $d_{j}:K\rightarrow \Xi(K/X)$ which we call the \emph{universal $j$-derivation}.
\end{definition}

\begin{definition}
Let $(K,D)$ be a $j$-field, let $X\subseteq K$, and $a\in K$. We say that $a$ belongs to the \emph{$j$-closure of $X$}, denoted $a\in j\mathrm{cl}(X)$, if for every $\partial\in j\mathrm{Der}(K/X)$ we have that $\partial(a)=0$. That is:
\begin{equation*}
    j\mathrm{cl}(X) = \bigcap_{\partial\in j\mathrm{Der}(K/X)}\ker\partial.
\end{equation*}
Combining the results of \cite[Lemmas 5.5, 5.6 and 5.7]{Eterovic-Schan-for-j} we know that $j\mathrm{cl}$ is a pregeometry. Let $\dim^{j}$ denote the dimension defined by $j\mathrm{cl}$. 
\end{definition}

The following proposition is a consequence of Ax-Schanuel for $j$ (Theorem \ref{j-chapter-Ax-for-j}).  
\begin{proposition}[see {{\cite[Proposition 6.2]{Eterovic-Schan-for-j}}}]
\label{prop:as}
Let $(K,D)$ be a $j$-field. Let $z_{1},\ldots,z_{n}\in D$ and let $F\subseteq K$ be $j\mathrm{cl}$-closed. Then:
\begin{equation*}
\mathrm{tr.deg.}_{F}F(\mathbf{z},J(\mathbf{z})) \geq 3\dim_{G}(\mathbf{z}|F) + \dim^{j}(\mathbf{z}|F).
\end{equation*}
\end{proposition}

We introduce a new piece of notation: we will sometimes use $j^{(t)}$ with $t=0,1,2,3$ to denote the functions $j$, $j'$, $j''$, and $j'''$ respectively. 

\begin{proposition}
\label{prop:c}
Let $(K,D)$ be a $j$-field, let $F\subseteq K$ be $j\mathrm{cl}$-closed, and let $z\in D$. 
\begin{enumerate}[(a)]
    \item If $z\in F$, then $\left\{j(z),j'(z),j''(z)\right\}\subseteq F$
    \item If $j^{(t)}(z)\in F$ for some $t\in\left\{0,1,2\right\}$, then $z\in F$.  
\end{enumerate}
\end{proposition}
\begin{proof}
By Proposition \ref{prop:as} we know that
 $$\mathrm{tr.deg.}_{F}F(z,J(z)) \geq 3\dim_{G}(z|F) + \dim^{j}(z|F).$$
 As $F$ is $j\mathrm{cl}$-closed, then the number on the right hand side of the inequality is either 0 or 4, depending on whether $z\in F$ or not. This completes the proof.
\end{proof}

\section{Finitely Generated \texorpdfstring{$j$}{j}-subfields}
\label{sec:fingenjfields}
Let $(K,D)$ be a $j$-field and let $(C,D_{C}) := (j\mathrm{cl}(\emptyset), D\cap j\mathrm{cl}(\emptyset))$. Let $\mathscr{K}_{j}$ denote the collection of $j$-subfields of $(K,D)$ which are finitely generated over $(C,D_{C})$. Throughout this section we fix a $j$-field $(K,D)$ and the corresponding collection $\mathscr{K}_{j}$.

We begin by defining a predimension function $\delta_{j}$, and then we define two standard objects known as ``self-sufficient extensions'' and ``self-sufficient closure''. We have not included all the proofs as they are mostly just repeating arguments that are already detailed in \cite{Aslanyan-adequate-predim} and \cite{Kirby-semiab} with appropriate name-changing; instead we give specific references in each case.

\subsection{Predimension}
\label{subsec:predim}
\begin{definition}
Suppose $(A,D_{A})$ and $(B,D_{B})$ are $j$-subfields of $(K,D)$ such that $(B,D_{B})$ is finitely generated over $(A,D_{A})$. In this case we define the \emph{$j$-predimension}:
$$\delta_{j}((B,D_{B})|(A,D_{A})) := \mathrm{tr.deg.}_{A}B - 3\dim_{G}(D_{B}|D_{A}).$$
More generally, given a subset $X\subset B$ we define $\delta_{j}(X|(A,D_{A})) := \delta_{j}(\left<X|(A,D_{A})\right>_{j}|(A,D_{A}))$.
When $(B,D_{B})\in\mathscr{K}_{j}$, we also define $\delta_{j}((B,D_{B})):= \delta_{j}((B,D_{B})|(C,D_{C}))$. 
\end{definition}

In particular, if $\mathbf{z}$ is a tuple from $D$ and $(A,D_{A})$ is a $j$-subfield of $(K,D)$, then 
$$\delta_{j}(\mathbf{z}|(A,D_{A})) = \mathrm{tr.deg.}_{A}A(\mathbf{z},J(\mathbf{z})) - 3\dim_{G}(\mathbf{z}|D_{A}).$$

We remark that $\delta_{j}$ is well-defined because, by definition, if  $(B,D_{B})$ is finitely generated over $(A,D_{A})$, then $\mathrm{tr.deg.}_{A}B$ and $\dim_{G}(D_{B}|D_{A})$ are both finite. 

\begin{lemma}
\label{lem:deltasubmodular}
The predimension $\delta_{j}$ is submodular, that is, for any $(A,D_{A}),(B,D_{B})$ in $\mathscr{K}_{j}$ we have:
$$\delta_{j}((A,D_{A})\cup(B,D_{B})) + \delta_{j}((A,D_{A})\cap(B,D_{B})) \leq \delta_{j}((A,D_{A})) + \delta_{j}((B,D_{B})).$$ 
\end{lemma}
\begin{proof}
On one hand it is easy to see that:
$$\dim_{G}(D_{A}\cup D_{B}|D_{C}) + \dim_{G}(D_{A}\cap D_{B}|D_{C}) = \dim_{G}(D_{A}|D_{C}) + \dim_{G}(D_{B}|D_{C}).$$
Similarly, it is a general property of field extensions that:
$$\mathrm{tr.deg.}_{C}(AB) + \mathrm{tr.deg.}_{C}(A\cap B)\leq \mathrm{tr.deg.}_{C}A + \mathrm{tr.deg.}_{C}B,$$
where $AB$ denotes the subfield of $K$ generated by $A\cup B$. 
\end{proof}

The main result of this section (Theorem \ref{thm:predim}) gives a characterisation of $\dim^{j}$ in terms of $\delta_{j}$. 

\subsection{Self-sufficient extensions}
\label{susbec:selfsufext}
We preserve the notation $(K,D)$, $(C,D_{C})$, and $\mathscr{K}_{j}$ from the beginning of the section. 

\begin{definition}
An embedding of $j$-fields $f:(K_{1},D_{1})\hookrightarrow (K_{2},D_{2})$ is called \emph{self-sufficient} (also known as \emph{strong}) if for every tuple $\mathbf{z}$ of $D_{2}$ we have that $\delta_{j}(\mathbf{z}|(K_{1},D_{1}))\geq 0$. We denote this property as $(K_{1},D_{1})\lhd  (K_{2},D_{2})$. We also say that $(K_{2},D_{2})$ is a \emph{self-sufficient extension} of $(K_{1},D_{1})$. 
\end{definition}

\begin{example}
\begin{enumerate}[(a)]
    \item The identity $\mathrm{id}:(K,D)\rightarrow(K,D)$ is a strong embedding. 
    \item Let $\mathbf{z}$ be a tuple of elements from $D$. By Proposition \ref{prop:as} we get that:
    $$\delta_{j}(\mathbf{z}|(C,D_{C}))\geq \dim^{j}(\mathbf{z}|D_{C})\geq 0.$$
    Therefore, the inclusion map $(C,D_{C})\hookrightarrow (K,D)$ is self-sufficient, or in symbols: $(C,D_{C})\lhd  (K,D)$. 
\end{enumerate}
\end{example}

\begin{lemma}
\label{lem:selfsuffequiv}
Let $(K_{1},D_{1})$ be a $j$-subfield of $(K,D)$ which contains $(C,D_{C})$. Then $(K_{1},D_{1})\lhd  (K,D)$ if and only if for every $(A,D_{A})\in\mathscr{K}_{j}$ we have $\delta_{j}((K_{1},D_{1})\cap (A,D_{A}))\leq \delta_{j}((A,D_{A}))$. 
\end{lemma}
\begin{proof}
Repeat proof of \cite[Lemma 2.9]{Aslanyan-adequate-predim}
\end{proof}

\begin{lemma}
\label{lem:selfsufftrans}
Let $(K_{1},D_{1}), (K_{2},D_{2}), (K_{3},D_{3})\in\mathscr{K}_{j}$. If $(K_{1},D_{1})\lhd  (K_{2},D_{2})$ and $(K_{2},D_{2})\lhd  (K_{3},D_{3})$, then $(K_{1},D_{1})\lhd  (K_{3},D_{3})$ (composition of self-sufficient embeddings of $j$-fields is self-sufficient).
\end{lemma}
\begin{proof}
Let $\mathbf{z}$ be a tuple of $D_{3}$. Observe that 
\begin{eqnarray*}
      \dim_{G}(\mathbf{z}|D_{1}) &=& \dim_{G}(\mathbf{z}|D_{2}) + \dim_{G}(\mathbf{z}\cap D_{2}|D_{1})\\
      \mathrm{tr.deg.}_{K_{1}}K_{1}(\mathbf{z},J(\mathbf{z})) &=& \mathrm{tr.deg.}_{K_{2}}K_{2}(\mathbf{z},J(\mathbf{z})) + \mathrm{tr.deg.}_{K_{1}}K_{1}(K_{2}\cap K_{1}(\mathbf{z},J(\mathbf{z}))).
\end{eqnarray*}
Let $\mathbf{z}'$ be the tuple of elements in $\mathbf{z}$ which are in $D_{2}$. Then the above equations show that:
$$\delta_{j}(\mathbf{z}|(K_{1},D_{1})) = \delta_{j}(\mathbf{z}|(K_{2},D_{2})) + \delta_{j}(\mathbf{z}'|(K_{1},D_{1}))\geq 0.$$
\end{proof}

\begin{definition}
Let $\lambda$ be an ordinal. A \emph{$\lambda$-chain of self-sufficient $j$-extensions} consists of a pair of families $\left\{(K_{\theta},D_{\theta})_{\theta<\lambda}, (f_{\theta_{1},\theta_{2}})_{\theta_{1}\leq\theta_{2}<\lambda}\right\}$ satisfying the following conditions:
\begin{enumerate}[(a)]
    \item For each $\theta<\lambda$, $(K_{\theta},D_{\theta})$ is a $j$-field.
    \item For each $\theta_{1}\leq\theta_{2}<\lambda$ the map $f_{\theta_{1},\theta_{2}}:(K_{\theta_{1}},D_{\theta_{1}})\rightarrow (K_{\theta_{2}},D_{\theta_{2}})$ is a self-sufficient embedding.
    \item For all $\theta_{1}\leq\theta_{2}\leq\theta_{3}<\lambda$ we have that $f_{\theta_{2},\theta_{3}}\circ f_{\theta_{1},\theta_{2}} = f_{\theta_{1},\theta_{3}}$. 
    \item For each $\theta <\lambda$, $f_{\theta,\theta}$ is the identity on $(K_{\theta},D_{\theta})$. 
\end{enumerate}
\end{definition}

\begin{lemma}
\label{lem:chainext}
Let $\lambda$ be an ordinal, let $\left\{(K_{\theta},D_{\theta})_{\theta<\lambda}, (f_{\theta_{1},\theta_{2}})_{\theta_{1}\leq\theta_{2}<\lambda}\right\}$ be a $\lambda$-chain of self-sufficient $j$-extensions, and let $(K,D)$ be the union of the chain.
\begin{enumerate}[(a)]
    \item Then $(K_{\theta},D_{\theta})\lhd  (K,D)$ for each $\theta<\lambda$.
    \item Suppose $(S,D_{S})$ is a $j$-field and that $(K_{\theta},D_{\theta})\lhd  (S,D_{S})$ for each $\theta<\lambda$. Then $(K,D)\lhd  (S,D_{S})$. 
\end{enumerate}
\end{lemma}
\begin{proof}
 Clear by finiteness arguments.
\end{proof}

\begin{proposition}[cf. {{\cite[Proposition 5.6]{Kirby-Schanuel}}}]
\label{prop:chain}
Suppose $(A,D_{A})\lhd (B,D_{B})$ is a self-sufficient inclusion of $j$-fields, both of which are graph-generated by $J$. Then there is an ordinal $\lambda$ and a $(\lambda+1)$-chain $\left\{(K_{\theta},D_{\theta})_{\theta\leq\lambda}, (f_{\theta_{1},\theta_{2}})_{\theta_{1}\leq\theta_{2}<\lambda}\right\}$ of self-sufficient $j$-extensions such that for all $0\leq\theta_{1}\leq\theta_{2}\leq\lambda$ we have:
\begin{enumerate}[(a)]
    \item $(A,D_{A}) = (K_{0},D_{0})$ and $(B,D_{B}) = (K_{\lambda},D_{\lambda})$,
    \item $K_{\theta_{1}}$ is graph-generated by $J$,
    \item For limit $\theta_{2}$, $(K_{\theta_{2}}, D_{\theta_{2}}) = \bigcup_{\theta_{1}<\theta_{2}}(K_{\theta_{1}},D_{\theta_{1}})$,
    \item $\dim_{G}(D_{\theta_{1}+1}|D_{\theta_{1}})$ and $\mathrm{tr.deg.}(K_{\theta_{1}+1}|K_{\theta_{1}})$ are finite,
    \item $(K_{\theta_{1}}, D_{\theta_{1}})\lhd (K_{\theta_{2}}, D_{\theta_{2}})$.
\end{enumerate}
\end{proposition}
\begin{proof}
  The proof is the same as that of \cite[Proposition 5.6]{Kirby-Schanuel} with the appropriate name-changing. 
\end{proof}

\subsection{Self-sufficient closure}

We now present a standard construction called the self-sufficient closure and show a few of its properties. For reference, see \cite[\textsection 2.1]{Aslanyan-adequate-predim} and \cite[\textsection 2.4]{Kirby-semiab}. We maintain the notation used at the beginning of \textsection \ref{sec:fingenjfields}. 

\begin{lemma}
\label{lem:selfsuffinter}
Let $\left\{(K_{i},D_{i})\right\}_{i}$ be a collection of $j$-subfields of $(K,D)$ all of which contain $(C,D_{C})$, and such that $(K_{i},D_{i})\lhd  (K,D)$ for all $i\in I$. Then 
$$\bigcap_{i\in I}(K_{i},D_{i})\lhd  (K,D).$$
\end{lemma}
\begin{proof}
Repeat the proof of \cite[Lemma 2.12]{Kirby-semiab}. 
\end{proof}

\begin{definition}
If $X$ is a subset of $K$, then the \emph{self-sufficient closure} of $X$ is the intersections of all $j$-subfields $(A,D_{A})$ of $(K,D)$ containing $\left<X\right>_{j}$ such that $(A,D_{A})\lhd  (K,D)$. We denote the self-sufficient closure of $X$ by $\lceil X \rceil$. 
\end{definition}

\begin{lemma}
\label{lem:selfsuffcl}
Let $(A,D_{A})\in\mathscr{K}_{j}$. Then:
\begin{enumerate}[(a)]
    \item $\lceil (A,D_{A}) \rceil\in\mathscr{K}_{j}$.
    \item $\lceil (A,D_{A}) \rceil\lhd  (K,D)$.
    \item $\delta_{j}\left(\lceil (A,D_{A}) \rceil\right) = \min\left\{\delta_{j}(B,D_{B}) : (A,D_{A})\subset(B,D_{B})\in\mathscr{K}_{j}\right\}$.
\end{enumerate}
\end{lemma}
\begin{proof}
For parts (a) and (c) repeat the arguments in \cite[Lemma 2.14]{Aslanyan-adequate-predim}. Part (b) follows immediately from Lemma \ref{lem:selfsuffinter}.
\end{proof}

\subsection{Extending \texorpdfstring{$j$}{j}-derivations}
\label{subsec:extender}
In this section we prove Theorem \ref{thm:predim}. For this, we first need to show that $j$-derivations can be extended in self-sufficient extensions; this is Proposition \ref{prop:extender2}. 

We point out that the proof of Proposition \ref{prop:extender2} crucially differs from its exponential counterpart (\cite[Theorem 6.3]{Kirby-Schanuel}) in the proof of Claim \ref{claim}. The proof of \cite[Theorem 6.3]{Kirby-Schanuel} uses intermediate results of Ax's proof of \cite[Theorem 3]{Ax} (the Ax-Schanuel theorem), and as Ax's proof is done with differential algebra, these intermediate steps are still valid in exponential fields. The proof of Ax-Schanuel for $j$ on the other hand (\cite{Pila-Tsim-Ax-j}) is done with o-minimality, and general $j$-fields have no o-minimal structure. 

Instead, for proving Proposition \ref{prop:extender2}, we will adapt an argument we presented in \cite[Theorem 3.5]{Aslanyan-Eterovic-Kirby-Diff-EC-j} as part of our solution of the differential version of EC for $j$, which just uses the statement of the Ax-Schanuel theorem. Furthermore, our method is rather general and can be expected to work for other functions for which there is a corresponding Ax-Schanuel theorem. For example, one can easily adapt our proof of Claim \ref{claim} to give a new proof of \cite[Theorem 6.3]{Kirby-Schanuel}.

\begin{proposition}
\label{prop:extender2}
Suppose $(A,D_{A})\lhd  (B,D_{B})$ is a self-sufficient extension of $j$-fields and that $A$ is graph-generated by $J$. Then every $j$-derivation on $(A,D_{A})$ extends to a $j$-derivation on $(B,D_{B})$. 
\end{proposition}
\begin{proof}
The proof is obtained from the proof of \cite[Theorem 3.5]{Aslanyan-Eterovic-Kirby-Diff-EC-j} after some appropriate reinterpretations. We give the full prof here as this is the main step in all of our main results.

  Let $K'$ be the subfield of $B$ generated by $D_{B}\cup J(D_{B})$, then $(K',D_{B})$ is a $j$-subfield of $(B,D_{B})$ which is graph-generated by $J$. Every $j$-derivation on $K'$ can be extended to $B$ by standard results of derivations, as this extension need only respect field operations. So we will assume that $(B,D_{B})$ is graph-generated by $J$. By Proposition \ref{prop:chain} we can further assume that $\dim_{G}(D_{B}|D_{A})$ and $\mathrm{tr.deg.}_{A}B$ are finite. Let $\partial_{0}$ be a $j$-derivation on $(A,D_{A})$. If $\partial_{0}=0$, then the result is trivial, so we assume that $\partial_{0}\neq 0$. Let $C_{0} = \ker\partial\subseteq A$. 
  
  Consider the space: 
  $$\mathrm{Der}(B|\partial_{0}):= \left\{\partial\in\mathrm{Der}(B|C_{0}) : \partial|_{A} = \lambda\partial_{0} \mbox{ for some } \lambda \in B\right\}.$$
  For every $\partial \in \Der(B|\partial_{0})$ there is a unique $\lambda_{\partial}\in B$ such that $\partial|_A = \lambda_{\partial} \partial_{0}$. The function $\varphi: \partial \mapsto \lambda_{\partial}$ gives a linear map $\varphi: \Der(B|\partial_{0})\rightarrow B$. This map is surjective, since for every $\lambda \in B$ the map $\lambda \partial_{0}:A\rightarrow B$ can be extended to a derivation of $B$. Moreover, $\ker(\varphi) = \Der(B|A)$, hence $$\dim \Der(B|\partial_{0}) = \dim \Der(B|A) +1 = \mathrm{tr.deg.}_{A}B+1.$$
  
  Consider the sequence of inclusions
$$\Der(B|A) \hookrightarrow \Der(B|\partial_{0}) \hookrightarrow \Der(B|C_{0}).$$
We then get a sequence of surjections
$$\Omega(B|C_{0}) \twoheadrightarrow \Omega(B|\partial_{0}) \twoheadrightarrow \Omega(B|A),$$ 
where $\Omega(B|\partial_{0})$ is defined as the dual of $\Der(B|\partial_{0})$. 
  
  Let $z_{1},\ldots,z_{n}$ be a $G\mathrm{cl}$-basis for $D_{B}$ over $D_{A}$. As $(A,D_{A})\lhd  (B,D_{B})$ we have that: 
  $$\mathrm{tr.deg.}_{A}B = \mathrm{tr.deg.}_{A}A(\mathbf{z},J(\mathbf{z}))\geq \dim_{G}(\mathbf{z}|D_{A}) = 3n.$$
  Let $\ell\geq 0$ be an integer such that $\dim \Omega(B|A) = 3n +\ell$. For $i\in\left\{1,\ldots,n\right\}$ define:
  \begin{equation*}
      \beta_{i} := d(j(z_{i})) - j'(z_{i})d(z_{i}),\quad \beta'_{i} := d(j'(z_{i})) - j''(z_{i})d(z_{i}),\quad \beta''_{i} := d(j''(z_{i})) - j'''(z_{i})d(z_{i}),
  \end{equation*}
  where $d:B\rightarrow\Omega(B|C_{0})$ is the universal derivation. Let $\Lambda(B|C_{0})$ be the $B$-linear subspace of $\Omega(B|C_{0})$ generated by $\left\{\beta_{i}, \beta'_{i}, \beta''_{i}\right\}_{i=1}^{n}$, and let $\Lambda(B|\partial_{0})\subseteq\Omega(B|\partial_{0})$ and $\Lambda(B|A)\subseteq\Omega(B|A)$ be the images of $\Lambda(B|C_{0})$ under the surjections above. 
  
\begin{claim}
\label{claim}
The forms $\beta_{i}, \beta_{i}', \beta_{i}'', i=1,\ldots, n$ are $B$-linearly independent in $\Omega(B|A)$, that is, $\dim \Lambda(B|A)  = 3n$.
\end{claim}
\begin{proof}
We proceed by contradiction, so assume $\dim \Lambda(B|A)  < 3n$. Consider the annihilator $\Ann(\Lambda(B|A)) \subseteq \Der(B|A)$, and observe that
\begin{equation}
    \label{eq:annhilator}
    \Ann(\Lambda(B|A)) = j\Der (B|A).
\end{equation} 
Clearly, 
$$r:=\dim \Ann(\Lambda(B|A)) = \dim \Omega(B|A) - \dim \Lambda(B|A) > \ell.$$
It is easy to see that $\Ann(\Lambda(B|A))$ is closed under the Lie bracket, hence we can choose a commuting basis of derivations $\partial_1,\ldots,\partial_r \in \Ann(\Lambda(B|A))$ (see \cite[Chapter 0, \S 5, Proposition 6]{Kolchin-diff-alg-gp} or \cite[Lemma 2.2]{Singer-noncommuting}). Let $L:= \bigcap_{i=1}^r\ker \partial_i$; thus $A \subseteq L \subsetneq B$. Also, by (\ref{eq:annhilator}), $L$ is the intersection of kernels of $j$-derivations on $B$, and so $L$ is $j\mathrm{cl}$-closed.

Let $\mathbf{v}_i:=(z_i,j(z_i),j'(z_i),j''(z_i))$. By Proposition \ref{prop:c}, either every coordinate of $\mathbf{v}_i$ is in $L$, or none of them are. Since $r> 0$, we may assume that for some $t\geq 1$ no coordinate of $\mathbf{v}_1, \ldots, \mathbf{v}_t$ is in $L$, and all coordinates of $\mathbf{v}_{t+1},\ldots, \mathbf{v}_n$ are in $L$. Let $$\mathbf{u}:= (\mathbf{v}_1, \ldots, \mathbf{v}_t),~ \mathbf{w} := (\mathbf{v}_{t+1},\ldots,\mathbf{v}_n).$$

As explained in the proof of the Claim of \cite[Theorem 3.5]{Aslanyan-Eterovic-Kirby-Diff-EC-j}, $\rk(\partial_i z_k)_{1\leq i\leq r,1\leq k\leq t}=r$. By the Ax-Schanuel theorem (Theorem \ref{j-chapter-Ax-for-j}) we get:
$$\mathrm{tr.deg.}_{L}L(\mathbf{u}) \geq 3t+\rk(\partial_i z_k)_{1\leq i\leq r,1\leq k\leq t} = 3t+r.$$ 
Further, using that $(A,D_{A})\lhd  (B,D_{B})$ we get
$$\mathrm{tr.deg.}_{A}L\geq \mathrm{tr.deg.}_{A}A(\mathbf{w})\geq 3\dim_{G}(z_{t+1},\ldots,z_{n}|D_{A}) = 3(n-t).$$ 
Combining these two inequalities we get 
$$\mathrm{tr.deg.}_{A}B = \mathrm{tr.deg.}_{L}B + \mathrm{tr.deg.}_{A}L \geq 3t+r + 3(n-t) = 3n+r > 3n + \ell,$$ 
which is a contradiction.
\end{proof}
  
  By Claim \ref{claim} the dimension of $\Lambda(B|\partial_{0})$ is also $3n$. Therefore, 
  $$\dim\mathrm{Ann}(\Lambda(B|\partial_{0})) = \dim \Omega(B|\partial_{0}) - \dim \Lambda(B|\partial_{0}) = 3n+\ell+1-3n =\ell+1$$ 
  and 
  $$\dim \Ann \Lambda(B|A) = \dim \Omega(B|A) - \dim \Lambda(B|A) = \ell.$$ 
  Choose a derivation $\partial \in \Ann \Lambda(B|\partial_{0}) \setminus \Ann \Lambda(B|A)$. Then $\partial|_A = \lambda_{\partial} \cdot \partial_{0}$ for some $\lambda_{\partial} \in B$. On the other hand,  $\partial \notin \Ann(\Lambda(B|A))$, therefore $\partial|_A \neq 0$ and $\lambda_{\partial} \neq 0$. Replacing $\partial$ by $\lambda_{\partial}^{-1}\cdot \partial$ we may assume that $\lambda_{\partial} = 1$ and $\partial$ is a $j$-derivation on $B$ which extends $\partial_{0}$. 
\end{proof}

The next result describes $\dim^{j}$ in terms of $\delta_{j}$. As is explained in \cite[\textsection 2.1]{Aslanyan-adequate-predim} and \cite[\textsection 2.7]{Kirby-semiab}, one can use the predimension $\delta_{j}$ to build a natural pregeometry and dimension on $K$. 

\begin{definition}
Given a finite set $X\subset K$ we define $\dim_{\delta_{j}}(X):= \delta_{j}\left(\lceil X\rceil\right)$, or equivalently 
$$\dim_{\delta_{j}}(X) = \min\left\{\delta_{j}(X\cup Y|(C,D_{C})) : Y \mbox{ is a finite subset of } K\right\}.$$
For $X$ as above and any subset $A\subset K$, we define:
$$\dim_{\delta_{j}}(X|A):=\min\left\{\dim_{\delta_{j}}(X\cup Y) - \dim_{\delta_{j}}(Y) | Y\mbox{ is a finite subset of } A\right\}.$$
Now, given any subset $B\subset K$, we define:
$$\mathrm{cl}_{\delta_{j}}(B) := \left\{x\in K : \dim_{\delta_{j}}(x|B) = 0\right\}.$$
\end{definition}

It follows easily from \cite[\textsection 2.1]{Aslanyan-adequate-predim} or \cite[Proposition 2.25]{Kirby-semiab} that $\mathrm{cl}_{\delta_{j}}$ is a pregeometry on $K$ whose corresponding dimension is $\dim_{\delta_{j}}$. Among other things, Theorem \ref{thm:predim} says that in fact $\dim_{\delta_{j}}$ agrees with $\dim^{j}$, and so $\mathrm{cl}_{\delta_{j}}$ agrees with $j\mathrm{cl}$. 

\begin{theorem}
\label{thm:predim}
Let $(K,D)$ be a $j$-field. For any tuple $\mathbf{x}$ of elements in $K$, we have that
$$\dim^{j}(\mathbf{x}) = \min\left\{\delta_{j}(\mathbf{x}\cup\mathbf{y}) : \mathbf{y} \mbox{ is a tuple from } K\right\}.$$
\end{theorem}
\begin{proof}
  Choose $\mathbf{y}$ such that $r:=\delta_{j}(\mathbf{x}\cup\mathbf{y})$ is minimal (we can do this by Proposition \ref{prop:as}). Let $(K_{0},D_{0}) = \left<\mathbf{x}\cup\mathbf{y}|(C,D_{C})\right>_{j}$. Let $\mathbf{z}$ be a $G\mathrm{cl}$-basis for $D_{0}$ over $D_{C}$. With the notation of the proof of Proposition \ref{prop:extender2}, by the observation in (\ref{eq:annhilator}) we have $j\mathrm{Der}(K_{0}|C) = \mathrm{Ann}(\Lambda(K_{0}|C))$, but by Claim \ref{claim}, $\mathrm{Ann}(\Lambda(K_{0}|C))$ has codimension $\dim_{G}(\mathbf{z}|D_{C})$ in $\mathrm{Der}(K_{0}|C)$. So:
  \begin{equation*}
      \dim j\mathrm{Der}(K_{0}|C) = \mathrm{tr.deg.}_{C}K_{0} - \dim_{G}(\mathbf{z}|D_{C}) = \delta_{j}(\mathbf{x}\cup\mathbf{y}) = r.
  \end{equation*}
  By the minimality of $r$ we get that, for every tuple $\mathbf{w}$ in $D$, 
  \begin{equation*}
      \delta_{j}(\mathbf{w}|(K_{0},D_{0})) = \delta_{j}(\mathbf{x}\cup\mathbf{w}\cup\mathbf{y})-\delta_{j}(\mathbf{x}\cup\mathbf{y})\geq 0, 
  \end{equation*}
  so $(K_{0},D_{0})\lhd  (K,D)$. By Proposition \ref{prop:extender2}, $j$-derivations on $K_{0}$ extend to $K$, which means that $\dim^{j}(\mathbf{x})\geq r$. By Proposition \ref{prop:as} we conclude then that $\dim^{j}(\mathbf{x}) = r$. 
\end{proof}

\section{Convenient Generators}
\label{sec:genJfinfields}
In \cite[Theorem 5.17]{Eterovic-Schan-for-j} it is shown that there exist many non-trivial $j$-derivations $\partial:\C\rightarrow\C$. Let $C = j\mathrm{cl}(\emptyset)\subset\mathbb{C}$. By the results of \cite[\textsection 5]{Eterovic-Schan-for-j}, $C$ is a countable algebraically closed subfield of $\mathbb{C}$, and if we set $D_C = \h\cap C$, then $(C,D_{C})$ is a full $j$-subfield of $(\C,\h)$.

\begin{theorem}
\label{thm:mainJ}
Let $F$ be a subfield of $\mathbb{C}$ such that $\mathrm{tr.deg.}_{C}F$ is finite. Then there exist  $t_{1},\ldots,t_{m}\in\mathbb{H}\setminus D_C$ such that:
\begin{enumerate}
    \item[(A1):] $F\subseteq\overline{C\left(\mathbf{t},J\left(\mathbf{t}\right)\right)}$, and
    \item[(A2):] $\mathrm{tr.deg.}_{C}C\left(\mathbf{t},J\left(\mathbf{t}\right)\right) = 3\dim_{G}\left(\mathbf{t}|C\right) + \dim^{j}\left(\mathbf{t}|C\right)$.
\end{enumerate} 
\end{theorem}
\begin{proof}
Let $T$ be a transcendence basis for $F$ over $C$. If $T$ is empty, then we are done as $F\subset C$. Otherwise, observe that $F$ is contained in a finite extension of $C(T)$. Let $(K,D)$ be the self-sufficient closure of $\left<j^{-1}(T)|(C,D_{C})\right>_{j}$ in $(\mathbb{C},\mathbb{H})$. Let $\mathbf{t}$ be a $G\mathrm{cl}$-basis for $D$ over $D_{C}$. Then we have that  
$$\mathrm{tr.deg.}_{C}C(\mathbf{t},J(\mathbf{t}))\leq\mathrm{tr.deg.}_{C}K.$$
Let $L$ be the relative algebraic closure of $C(\mathbf{t},J(\mathbf{t}))$ in $K$ and consider the $j$-field $(L,D)$. Then $\delta_{j}((L,D))\leq \delta_{j}((K,D))$. By Lemma \ref{lem:selfsuffcl} we know that $(K,D)$ minimizes the predimension $\delta_{j}$, so $\delta_{j}((L,D)) = \delta_{j}((K,D))$, and as $\dim_{G}(D|D_{C}) = \dim_{G}(\mathbf{t})$, we conclude that $\mathrm{tr.deg.}_{C}L = \mathrm{tr.deg.}_{C} K$. So in fact $L=K$. This verifies condition (A1). 

By Theorem \ref{thm:predim} and Lemma \ref{lem:selfsuffcl} we get that 
$$\dim^{j}(\mathbf{t}) = \dim^{j}(K|C)=\dim^{j}(K) = \delta_{j}((K,D)),$$ 
which verifies (A2). 
\end{proof}

In the proof of Theorem \ref{thm:mainJ}, the elements $t_{1},\ldots,t_{m}\in\mathbb{H}\setminus D_C$ come from the self-sufficient closure of $F$, in which case the following three results are immediate by Lemma \ref{lem:selfsuffcl}. But on its own, condition (A2) still allows us to get a version of Proposition \ref{prop:as} even if the field is not $j\mathrm{cl}$ closed. 

\begin{lemma}
\label{lem::asgJf}
Suppose $t_{1},\ldots,t_{m}\in\h\setminus D_{C}$ satisfy condition (A2). Then for any $z_{1},\ldots,z_{n}\in \mathbb{H}$ we have:
$$\mathrm{tr.deg.}_{C(\mathbf{t},J(\mathbf{t}))}C(\mathbf{z},\mathbf{t},J(\mathbf{z}),J(\mathbf{t}))\geq3\dim_{G}(\mathbf{z}|C\cup\mathbf{t}) + \dim^{j}(\mathbf{z}|C\cup\mathbf{t}).$$
\end{lemma}
\begin{proof}
The inequality is obtained by first using Proposition \ref{prop:as} to get:
$$\mathrm{tr.deg.}_{C}C(\mathbf{z},\mathbf{t},J(\mathbf{z}),J(\mathbf{t}))\geq3\dim_{G}(\mathbf{z}\cup\mathbf{t}|C) + \dim^{j}(\mathbf{z}\cup\mathbf{t}|C)$$
and now using the addition formula (see \cite[C.1.8]{tent-ziegler}) and (A2). 
\end{proof}

\begin{corollary}
Suppose $t_{1},\ldots,t_{m}\in\mathbb{H}\setminus D_C$ satisfy condition (A2). Let $(A,D_{A})$ be the $j$-subfield of $(\C,\h)$ generated by $D_{C}\cup\mathbf{t}$. Then $(A,D_{A})\lhd  (\C,\h)$.  
\end{corollary}
\begin{proof}
Given a tuple $\mathbf{z}$ of $\h$, Lemma \ref{lem::asgJf} implies that: $\delta_{j}(\mathbf{z}|(A,D_{A}))\geq \dim^{j}(\mathbf{z}|C\cup\mathbf{t})\geq 0$.
\end{proof}

We can also get a version of Lemma \ref{lem::asgJf} without derivatives. 

\begin{corollary}
\label{cor:genjfin}
Suppose $t_{1},\ldots,t_{m}\in\h\setminus D_{C}$ satisfy condition (A2). Then for any $z_{1},\ldots,z_{n}\in \mathbb{H}$ we have:
$$\mathrm{tr.deg.}_{C(\mathbf{t},J(\mathbf{t}))}C(\mathbf{z},\mathbf{t},j(\mathbf{z}),J(\mathbf{t}))\geq\dim_{G}(\mathbf{z}|C\cup\mathbf{t}) + \dim^{j}(\mathbf{z}|C\cup\mathbf{t}).$$
\end{corollary}

\subsection{Proof of Theorem \ref{thm:unconditional}}
We recall that \cite[Theorem 1.2]{Eterovic-Herrero} shows that if one assumes the modular version of Schanuel's conjecture (see \textsection \ref{subsec:msc}), then in every plane irreducible curve $V\subset\C^{2}$ which is not a horizontal or vertical line, there exists a generic point which is of the form $(z,j(z))$ for some $z\in\h$. Theorem \ref{thm:unconditional} gives unconditional cases of \cite[Theorem 1.2]{Eterovic-Herrero}. 

\begin{proof}[Proof of Theorem \ref{thm:unconditional}]
By \cite[Theorem 1.1]{Eterovic-Herrero} we know that there exist infinitely many $z\in\h$ such that $(z,j(z))\in V$. In particular, the set $S:=\left\{(z,j(z))\in V : z\in\h\right\}$ is Zariski dense in $V$. 

Let $F\subset\C$ be a finitely generated field over which $V$ is defined. By Theorem \ref{thm:mainJ} we know that there exist $t_{1},\ldots,t_{m}\in\h\setminus D_{C}$ such that by Corollary \ref{cor:genjfin} we have 
$$\mathrm{tr.deg.}_{F}F(z,j(z))\geq\dim_{G}(z|C\cup\mathbf{t}) + \dim^{j}(z|C\cup\mathbf{t})$$
for every $z\in\h$.

Now, by \cite[Proposition 7.13]{Eterovic-Herrero} there can only be finitely many elements $(z,j(z))\in S$ such that $z\in G\mathrm{cl}(\mathbf{t})$. Furthermore, as $S$ is Zariski dense in $V$ and $V$ is not defined over $C$, then by Proposition \ref{prop:c} we have that there are only finitely many elements in $S$ with coordinates in $C$. Therefore there exists $z\in\h$ such that $z\notin D_{C}\cup\mathbf{t}$ and $(z,j(z))\in V$. Therefore, by Corollary \ref{cor:genjfin}
$$1\leq \dim_{G}(z|C\cup\mathbf{t}) + \dim^{j}(z|C\cup\mathbf{t})\leq \mathrm{tr.deg.}_{F}F(z,j(z)) \leq \dim V = 1.$$
\end{proof}

What happens if $V$ is defined over $C$? As the next lemma shows, in that case any point $(z,j(z))\in V$ will have coordinates in $C$, and so Lemma \ref{lem::asgJf} can only provide trivial inequalities. We will show in \textsection\ref{subsec:Jfinfields} that if one assumes a modular version of Schanuel's conjecture, then we can get results like Theorem \ref{thm:mainJ} for subfield of $C$.

\begin{lemma}
Let $F$ be a subfield of $\C$ and let $p(X,Y)\in F[X,Y]$ be irreducible in $\mathbb{C}[X,Y]$. Suppose that $z\in\h$ is such that $p(z,j(z))=0$ and $\mathrm{tr.deg.}_{F}F(z,j(z))=1$. Then $z,j(z)\in j\mathrm{cl}(F)$. 
\end{lemma}
\begin{proof}
Since $p(z,j(z))=0$, then $\mathrm{tr.deg.}_{F}F(z,j(z),j'(z),j''(z))$ is at most 3. On the other hand, by Proposition \ref{prop:as} we have that $\mathrm{tr.deg.}_{j\mathrm{cl}(F)}j\mathrm{cl}(F)(z,j(z),j'(z),j''(z))$ is either 0 or 4, so it must by 0.
\end{proof}

\subsection{The exponential case}

We can adapt the ideas we have presented in this section to obtain unconditional results for the complex exponential function $\exp$ regarding the strong exponential closedness conjecture. An \emph{exponential derivation} on $\mathbb{C}$ is a field derivation $\partial:\mathbb{C}\rightarrow\mathbb{C}$ satisfying that $\partial(\exp(z)) = \exp(z)\partial(z)$ for every $z\in\mathbb{C}$. The fact that exponential derivations define a pregeometry $\mathrm{ecl}$ on every exponential field and the fact that it agrees with the pregeometry coming from a corresponding predimension, is given by the main results of \cite{Kirby-Schanuel} (the Ax-Schanuel theorem for the exponential function is given in \cite{Ax}). Let $E=\mathrm{ecl}(\emptyset)$, let $\dim^{e}$ denote the dimension coming from $\mathrm{ecl}$, and given $A,B\subset\mathbb{C}$, let $\mathrm{lin.dim}_{\mathbb{Q}}(A|B)$ denote the linear dimension of the $\mathbb{Q}$-vector space generated by $A\cup B$ modulo the subvector space generated by $B$. Repeating the proof of Theorem \ref{thm:mainJ}, we get the corresponding result on the existence of convenient generators (necessary results regarding self-sufficient closures in exponential fields can be found in \cite{Kirby-semiab}).
\begin{theorem}
\label{thm:maine}
Let $F$ be a subfield of $\mathbb{C}$ such that $\mathrm{tr.deg.}_{E}F$ is finite. Then there exist $t_{1},\ldots,t_{m}\in\C$ such that:
\begin{enumerate}
    \item[(E1):] $F\subseteq\overline{E\left(\mathbf{t},\exp\left(\mathbf{t}\right)\right)}$,
    \item[(E2):] $\mathrm{tr.deg.}_{E}E\left(\mathbf{t},\exp\left(\mathbf{t}\right)\right) = \mathrm{lin.dim}_{\mathbb{Q}}\left(\mathbf{t}|E\right) + \dim^{e}\left(\mathbf{t}|E\right)$.
\end{enumerate}
\end{theorem}

Proceeding just as we did in the proof of Theorem \ref{thm:unconditional} we get:
\begin{theorem} \label{thm:unconditionalfore}
Let $V\subset\C^{2}$ be an algebraic curve which is neither a horizontal nor a vertical line. If $V$ is not definable over $E$, then for any finitely generated field $F$ over which $V$ is defined, there exists $z\in\C$ such that $(z,\exp(z))\in V$ and $\mathrm{tr.deg.}_{F}F(z,\exp(z)) = 1$. 
\end{theorem}
\begin{proof}[Proof sketch]
By Theorem \ref{thm:maine}, we know that $F$ is contained in a field of the form $E\left(\mathbf{t},\exp\left(\mathbf{t}\right)\right)$ satisfying conditions (E1) and (E2). From \cite[Theorem 1.3]{mantova} we get that there are only finitely many vectors $\mathbf{c}\in\mathbb{Q}^{\ell}$ such that $(\mathbf{c}\cdot\mathbf{t},\exp(\mathbf{c}\cdot\mathbf{t}))\in V$. On the other hand, $V$ has a Zariski dense set of points of the form $(z,\exp(z))$ (this is obtained by using Hadamard's factorisation theorem, as explained in \cite{marker}). 
\end{proof}

\section{\texorpdfstring{$j$}{j}-polynomials}
\label{sec:jpolys}
In this section we show that, using our result on extension of $j$-derivations (Proposition \ref{prop:extender2}), we can give yet another characterisation of $j\mathrm{cl}$, this time in terms of systems of equations called \emph{Khovanskii systems}. More precisely, we will show that given $a\in K$ and $A\subseteq K$, then $a\in j\mathrm{cl}(A)$ if and only if $a$ satisfies a certain system of equations with coefficients coming from $A$ (see Theorem \ref{thm:extender3}). This characterisation will then be used to show  Theorem \ref{thm:mainJ2}, which says that, under a modular version of Schanuel's conjecture, we can find convenient sets of generators for some fields, in analogy with Theorem \ref{thm:mainJ}. 

\subsection{Khovanskii systems}
\label{subsec:khovanskii}
Khovanskii systems on $j$-fields are a straightforward analogue of \cite[\textsection 3 and \textsection 4]{Kirby-Schanuel}. The results are restated here in terms of $j$-fields; the proofs are mostly the same as in \cite{Kirby-Schanuel} with appropriate substitutions. 

\begin{definition}
Given a ring $R$ (commutative and unital), we define the ring $R[\mathbf{X}]^{J}$ of \emph{$j$-polynomials} on the variables $\mathbf{X} = (X_{1},\ldots,X_{n})$ as the ring
\begin{equation*}
    R[\mathbf{X}]^{J}:= R[\mathbf{X},j(\mathbf{X}),j'(\mathbf{X}),j''(\mathbf{X})] = R[\mathbf{X},J(\mathbf{X})].
\end{equation*}
 Note that the elements of this ring are formal expressions, but when we work in a $j$-field $(K,D)$ and consider the ring $K[\mathbf{X}]^{J}$, we will want to evaluate these expressions at some tuple $\mathbf{a}$ of $K$. This can only be done if the appropriate coordinates of $\mathbf{a}$ are in $D$. 
 
Given $i\in\left\{1,\ldots,n\right\}$, we define the operator $\frac{\partial}{\partial X_{i}}:R[\mathbf{X}]^{J}\rightarrow R[\mathbf{X}]^{J}[j'''(\mathbf{X})]$ as expected:\footnote{From the differential equation $\Psi(j,j',j'',j''')=0$, one sees that the ring $R[\mathbf{X}]^{J}[j'''(\mathbf{X})]$ is a subring of the fraction field of $R[\mathbf{X}]^{J}$.}
 \begin{enumerate}[(a)]
     \item For every $k\in \left\{1,\ldots,n\right\}$ we have
     $$\frac{\partial X_{k}}{\partial X_{i}} =\left\{\begin{array}{cc}
         0 &  \mbox{if } k\neq i\\
         1 & \mbox{if } k=i
     \end{array}\right. .$$
     \item For all $a,b\in R$ and all $f_{1},f_{2}\in R[\mathbf{X}]^{J}$:
     $$\frac{\partial (af_{1}(\mathbf{X}) + bf_{2}(\mathbf{X}))}{\partial X_{i}} = a\frac{\partial f_{1}(\mathbf{X})}{\partial X_{i}} + b\frac{\partial f_{2}(\mathbf{X})}{\partial X_{i}}.$$
     \item For all  $f_{1},f_{2}\in R[\mathbf{X}]^{J}$:
     $$\frac{\partial \left(f_{1}(\mathbf{X})f_{2}(\mathbf{X})\right)}{\partial X_{i}} = f_{2}(\mathbf{X})\frac{\partial f_{1}(\mathbf{X})}{\partial X_{i}} + f_{1}(\mathbf{X})\frac{\partial f_{2}(\mathbf{X})}{\partial X_{i}}.$$
     \item For every $k\in \left\{1,\ldots,n\right\}$ we have 
     $$\frac{\partial j^{(t)}\left(X_{k}\right)}{\partial X_{i}} =\left\{\begin{array}{cc}
         0 & \mbox{if } k\neq i\\
         j^{(t+1)}\left(X_{i}\right) & \mbox{if } k=i
     \end{array}\right. ,$$
     where $t=0,1,2$.
 \end{enumerate}
\end{definition}

\begin{definition}
Let $(K,D)$ be a $j$-field. Given a subset $B\subset K$, let $R_{B}$ denote the subring of $K$ generated by $B$. A \emph{Khovanskii system of $j$-polynomials over $B$} (in the variables $X_{1},\ldots,X_{n}$) consists of a set of $j$-polynomials $f_{1},\ldots,f_{n}\in R_{B}[X_{1},\ldots,X_{n}]^{J}$ (for some $n\in\mathbb{N}$), the system of equations:
\begin{equation*}
    f_{i}(X_{1},\ldots,X_{n}) = 0,\, \mbox{ for } i =1,\ldots,n,
\end{equation*}
and the inequation
\begin{equation*}
    \det\left[
    \frac{\partial f_{i}(\mathbf{X})}{\partial X_{k}}\right]
    _{i,k = 1,\ldots,n}\neq 0.
\end{equation*}
\end{definition}

\begin{definition}
Let $(K,D)$ be a $j$-field and let $B\subset K$ be any subset. We say that $a\in K$ belongs to the \emph{$k$-closure} of $B$ (written $a\in \mathrm{kcl}(B)$) if for some $n\in\mathbb{N}$ there exist $a_{1},\ldots,a_{n}\in K$, with $a_{1}=a$, and $j$-polynomials $f_{1},\ldots,f_{n}\in R_{B}[X_{1},\ldots,X_{n}]^{J}$, such that $\mathbf{a}$ satisfies the Khovanskii system determined by  $f_{1},\ldots,f_{n}$. 
\end{definition}

The main result of this section (Theorem \ref{thm:extender3}) is that $\mathrm{kcl}$ and $j\mathrm{cl}$ are actually the same operator. 

\begin{remark}
Perhaps some readers may have noticed that in our definition of the ring $R[\mathbf{X}]^{J}$ we have not included any expressions that include iterations of the functions $j$, $j'$ and $j''$. For example, consider the equation: $j( j'(X^{2}) + 4) = 1$. Even if this is not an equation that can be written with elements of $R[\mathbf{X}]^{J}$, we can however find a system of equations with elements of $R[\mathbf{X}]^{J}$ that will have the same solutions:
\begin{eqnarray*}
j(X_{1}) &=& 1\\
j'(X_{2}) + 4 &=& X_{1}\\
X_{3}^{2} &=& X_{2}.
\end{eqnarray*}
\end{remark}

The following lemma is straightforward (cf. \cite[Lemma 3.3]{Kirby-Schanuel}).

\begin{lemma}
\label{lem:kclosure}
For any $A,B\subseteq K$, the operator $\mathrm{kcl}$ satisfies:
\begin{enumerate}[(a)]
    \item $A\subseteq \mathrm{kcl}(A)$.
    \item $A\subseteq B\implies \mathrm{kcl}(A)\subseteq \mathrm{kcl}(B)$.
    \item $\mathrm{kcl}(\mathrm{kcl}(A)) = \mathrm{kcl}(A)$.
    \item $\mathrm{kcl}$ has finite character. 
\end{enumerate}
\end{lemma}

\begin{remark}
\label{rem:kclsolinD}
Let $(K,D)$ be a full $j$-field, and suppose $a_{1},\ldots,a_{n}\in K$ form a solution of some Khovanskii system of $j$-polynomials. Then we can assume that $a_{1},\ldots,a_{n}\in D$. Indeed, if $a_{2}\notin D$ say, then as $(K,D)$ is full there exists $b_{2}\in D$ such that $j(b_{2}) = a_{2}$. So, if $f(\mathbf{X})$ is a $j$-polynomial such that $f(\mathbf{a}) =0$, the system 
\begin{equation*}
    \begin{array}{ccc}
        f(X_{1}, Y, X_{3}, \ldots,X_{n}) & = & 0 \\
       j(X_{2}) & = & Y 
    \end{array}
\end{equation*}
is a system of $j$-polynomials which vanishes on $(a_{1},b_{2},a_{3},\ldots,a_{n})$. 
\end{remark}

\begin{remark}
\label{rem:kclalgcl}
Note that if $a\in K$ is algebraic over a subset $A$, then $a\in \mathrm{kcl}(A)$. Indeed, the minimal polynomial of $a$ with coefficients in $R_{A}$ (the ring generated by $A$) satisfies the conditions of a Khovanskii system. Therefore, $\mathrm{kcl}(A)$ is a relatively algebraically closed subfield of $K$. 

It will be convenient to also point out that in the case of full $j$-fields we could have defined $j$-polynomials in a different way, without changing the resulting $\mathrm{kcl}$ operator. We could have said that a $j$-polynomial over $R$ is any element of the ring
\begin{equation*}
    R\left[\left\{g\mathbf{X},J(g\mathbf{X})\right\}_{g\in G}\right],
\end{equation*}
where we treat the expressions $gX_{i}$ and $J(gX_{i})$ as abstract symbols until we decide to evaluate them at some point. Of course, we then need to redefine the operators $\frac{\partial}{\partial X_{i}}$, but this can still be done in the natural way. The point now is that when working over a $j$-field $(K,D)$, for every $a\in D$ and every $g\in G$ we have that $ga$ is algebraic over $\mathbb{Q}(a)$ and the coordinates of $J(ga)$ are algebraic over $\mathbb{Q}(a,J(a))$. For this reason (using Remark \ref{rem:kclsolinD}), defining $\mathrm{kcl}$ in the way we did or using this other definition, would not affect the operator. 
\end{remark}

Recall that (Remark \ref{rem:cdder}), given a $j$-field $(K,D)$, for any subset $A\subset K$ we have a universal derivation $d:K\rightarrow\Omega(K/A)$, and a universal $j$-derivation $d_{j}:K\rightarrow\Xi(K/A)$ which vanishes on $A$, where $d_{j}$ is just the composition of $d$ and the quotient map $\Omega(K/A)\rightarrow\Xi(K/A)$. The space $\Xi(K/A)$ satisfies the universal property that for any $\partial\in j\mathrm{Der}(K/A, M)$ there exists a $K$-linear map $\partial^{\ast}$ giving a commutative diagram:
   \begin{center}
    \begin{tikzcd}[ampersand replacement=\&]
    K \arrow[d,"\partial"]\arrow[r,"d"] \& \Xi(K/A) \arrow[dl,"\partial^{\ast}"]\\
    M \& 
    \end{tikzcd} 
    \end{center}
   which allows us to identify $j\mathrm{Der}(K/A)$ with the dual of $\Xi(K/X)$, and so we get $\dim\Xi(K/F) = \dim j\mathrm{Der}(K/F)$. The next lemma  gives a more concrete description of $\Xi(K/A)$. 

\begin{lemma}
\label{lem:univcharac}
Let $(K,D)$ be a $j$-field and let $A\subset K$. Let $M$ be the $K$-vector space generated by the symbols  $\left\{mr : r\in K\right\}$ subject to the relations:
\begin{equation}
\label{eq:derivation}
    \sum_{i=1}^{n}\frac{\partial f}{\partial X_{i}}(\mathbf{r})mr_{i} =0 
\end{equation}
for each $f\in R_{A}[\mathbf{X}]^{J}$ and every tuple $\mathbf{r}$ from $K$ such that $f(\mathbf{r})=0$. Then there exists an isomorphism  of $K$-vector spaces $s:M\rightarrow\Xi(K/A)$ satisfying $mr\mapsto d_{j}r$ for every $r\in K$. 
\end{lemma}
\begin{proof}
Straightforward adaptation of \cite[Lemma 4.6]{Kirby-Schanuel}.
\end{proof}

\begin{proposition}
\label{prop:kleqj}
Let $(K,D)$ be a $j$-field and $A\subset K$. Then $\mathrm{kcl}(A)\subseteq j\mathrm{cl}(A)$. 
\end{proposition}
\begin{proof}
Straightforward adaptation of \cite[Proposition 4.7]{Kirby-Schanuel}.
\end{proof}

\begin{lemma}
\label{lem:dergenerator}
Let $(K,D)$ be a $j$-field and let $(F,D_{F})$ be a $j$-subfield. Suppose that $z_{1},\ldots,z_{n}\in D$ form a $G\mathrm{cl}$-basis for $D$ over $D_{F}$, and that $K$ is graph-generated by $j$. Let $M$ be the $K$-vector space generated by the symbols $mz_{1},\ldots, mz_{n}$ subject to the relations:
\begin{equation}
\label{eq:dergenerator}
    \sum_{i=1}^{n}\frac{\partial f}{\partial X_{i}}(\mathbf{z})mz_{i} =0 
\end{equation}
for each $f\in F[\mathbf{X}]^{J}$ satisfying $f(\mathbf{z})=0$. Then there exists an isomorphism of $K$-vector spaces $s:M\rightarrow\Xi(K/F)$ satisfying $mz_{i}\mapsto d_{j}z_{i}$ for all $i\in\left\{1,\ldots,n\right\}$. 
\end{lemma}
\begin{proof}
Straightforward adaptation of \cite[Lemma 4.8]{Kirby-Schanuel}.
\end{proof}

For the following lemma, let $A^{\mathrm{alg}_K}$ denote the relative algebraic closure of the field generated by $A$ in $K$, that is for $A\subset K$, $A^{\mathrm{alg}_K}$ consists of all the elements of $K$ which are algebraic over the field generated by $A$. 

\begin{lemma}
\label{lem:graphgen}
Let $(K,D)$ be a $j$-field and let $A\subseteq K$ be a subfield. Let $A_{1} = A\cap\mathbb{Q}(D\cup J(D))^{\mathrm{alg}_K}$ and $A_{2} = A\setminus A_{1}$. Then
$$j\mathrm{cl}(A) = (j\mathrm{cl}(A_{1})\cup A_{2})^{\mathrm{alg}_K}.$$
\end{lemma}
\begin{proof}
Notice that $j$-derivations behave just like normal field derivations when applied to elements outside of the graph of $J$, so any element of $j\mathrm{cl}(A)$ which is not in $j\mathrm{cl}(A_{1})$ must be algebraic over $j\mathrm{cl}(A_{1})\cup A_{2}$. 
\end{proof}

In the proof of the following theorem, we will relativise the closure operators $\mathrm{kcl}$ and $j\mathrm{cl}$ to $j$-subfields. So, if $(F,D_{F})$ is a $j$-subfield of $(K,D)$ and $A$ is a subset of $F$, we will use the notation $j\mathrm{cl}_{F}(A):=j\mathrm{cl}(A)\cap F$ and $\mathrm{kcl}_{F}(A):=\mathrm{kcl}(A)\cap F$.

\begin{theorem}
\label{thm:extender3}
Let $(K,D)$ be a $j$-field and let $A\subseteq K$. Then $j\mathrm{cl}(A) = \mathrm{kcl}(A)$. 
\end{theorem}
\begin{proof}
  By Proposition \ref{prop:kleqj} we already have that $\mathrm{kcl}(A)\subseteq j\mathrm{cl}(A)$, so now we focus on proving the reverse inclusion. 
  
  Let $a\in j\mathrm{cl}(A)$. First we make some reductions. Because both $j\mathrm{cl}(A)$ and $\mathrm{kcl}(A)$ are relatively algebraically closed in $K$, we can assume that $a$ is not algebraic over $\mathbb{Q}$, and so in particular, we assume that $a\notin\Sigma$. Furthermore, by the finite character of $j\mathrm{cl}$, we may also assume that $A$ is finite. Finally, by Lemma \ref{lem:graphgen} we can assume that $K$ is graph-generated by $J$, and so it suffices to prove the theorem when $A\subset D$ and $a\in D$. 
  
  Let $(F,D_{F})$ be a finitely generated $j$-subfield  of $(K,D)$ which is graph-generated by $J$, and such that $A\subseteq D_{F}$, $a\in D_{F}$ and $a\in j\mathrm{cl}_{F}(A)$ (for example $(F,D_{F})$ can be the self-sufficient closure of $A\cup\left\{a\right\}$, see also \cite[Lemma 5.6]{Eterovic-Schan-for-j}). Furthermore, we choose $(F,D_{F})$ so that $\dim_{G}(D_{F}|A)$ is minimal with these properties. Combined, these statements will imply that $F = j\mathrm{cl}_{F}(A)$. We explain this in the following paragraph.
  
  Let $L = j\mathrm{cl}_{F}(A)$ and $D_{L} = D_{F}\cap L$, this way $(L,D_{L})$ is a $j$-subfield of $(F,D_{F})$. If $F\neq L$, then by the minimality of $F$ we have that $a\notin j\mathrm{cl}_{L}(A)$. Therefore there is a $j$-derivation $\partial\in j\mathrm{Der}(L|A)$ which does not extend to $F$. By Proposition \ref{prop:extender2} we conclude that $(F,D_{F})$ is not a strong extension of $(L,D_{L})$. However, this contradicts Proposition \ref{prop:as} because, as $L$ is $j\mathrm{cl}_{F}$-closed in $F$, for every tuple $\mathbf{z}$ in $D_{F}$, we have that $\delta_{j}(\mathbf{z}|D_{L})\geq 0$. In conclusion, $F=L$. 
  
  Now let $a=z_{1},\ldots,z_{n}$ be a $G\mathrm{cl}$-basis for $D_{F}$ over $D_{L}$. By Lemma \ref{lem:dergenerator}, $\Xi(F/A)$ is generated by the elements $d_{j}z_{1},\ldots,d_{j}z_{n}$ subject to the relations
  \begin{equation*}
      \sum_{i=1}^{n}\frac{\partial f}{\partial X_{i}}(\mathbf{z})d_{j}z_{i} =0, 
  \end{equation*}
  for each $f\in F[\mathbf{X}]^{J}$ satisfying $f(\mathbf{z})=0$ in $F$. As $F = j\mathrm{cl}_{F}(A)$, then $\Xi(F/A)=0$. This means that we can choose $f_{1},\ldots,f_{n}\in F[\mathbf{X}]^{J}$ such that $f_{i}(\mathbf{z})=0$ for every $i\in\left\{1,\ldots,n\right\}$, and the matrix $J = \left(\frac{\partial f_{t}}{\partial X_{s}}(\mathbf{z})\right)$ is non-singular. Therefore $a\in \mathrm{kcl}_{F}(A)$. As $\mathrm{kcl}_{F}(A)\subseteq \mathrm{kcl}(A)$, we are done. 
\end{proof}

Combined with Theorem \ref{thm:predim}, we have now proven Theorem \ref{thm:main}. In particular, Theorem \ref{thm:extender3} proves that $\mathrm{kcl}$ is a pregeometry.

\subsection{Higher dimensional examples of EC}
\label{subsec:higherdimec}
Let $C=j\mathrm{cl}(\emptyset)$ and choose $a\in\mathbb{C}\setminus C$. Given a positive integer $n$, let $V\subset\mathbb{C}^{2n}$ be the algebraic variety defined by the equations:
\begin{equation*}
    V = \left\{\begin{array}{ccc}
         X_{1} &=& Y_{n} + a  \\
         X_{2} &=& Y_{1}\\
         X_{3} &=& Y_{2} \\
         &\vdots& \\
         X_{n} &=& Y_{n-1}
    \end{array}\right\}.
\end{equation*}
\begin{proposition}
\label{prop:unconditonalhigherdim}
For every finitely generated field $K$, there exists a point $(\mathbf{z},j(\mathbf{z}))\in V$ such that $\mathrm{tr.deg.}_{K}K(\mathbf{z},J(\mathbf{z}))\geq n$. 
\end{proposition}
\begin{proof}
 By \cite[Theorem 1.1]{Eterovic-Herrero}, we know that $V$ has a Zariski dense set of points of the form $(\mathbf{z},j(\mathbf{z}))$. Let us define inductively the functions $j_{n}$ by: $j_{1}(z):=j(z)$ and $j_{k+1}(z):=j_{k}(j(z))$ for every $k>1$ (cf \cite[\textsection 2]{Eterovic-Herrero}). Then every solution $(\mathbf{z},j(\mathbf{z}))\in V$ satisfies that $z_{1} = j_{n}(z_{1})+a$. Furthermore, as $j(\h\cap C) = C$, then by Proposition \ref{prop:c} we must have that $z,j_{1}(z),j_{2}(z),\ldots,j_{n}(z)\notin C$, as otherwise we would contradict that $a\notin C$. 

By Theorem \ref{thm:mainJ}, there exist $t_{1},\ldots,t_{m}\in\h\setminus C$ such that conditions (g$J$f1) and (g$J$f2) are satisfied with respect to $K$. By Corollary \ref{cor:genjfin} we get that for every $(\mathbf{z},j(\mathbf{z}))\in V$:
$$\mathrm{tr.deg.}_{K}K(\mathbf{z},j(\mathbf{z}))\geq\dim_{G}(\mathbf{z}|C\cup\mathbf{t}) + \dim ^{j}(\mathbf{z}|C\cup\mathbf{t}).$$
Using that the points of the form $(\mathbf{z},j(\mathbf{z}))$ are Zariski dense in $V$ and \cite[Proposition 7.13]{Eterovic-Herrero}, we get that $V$ has a Zariski dense set of points of the form $(\mathbf{z},j(\mathbf{z}))$ such that $z_{1},\ldots,z_{n}\notin C\cup G\cdot\mathbf{t}$. For these solutions we then have that $\dim_{G}(\mathbf{z}|C\cup\mathbf{t}) = \dim_{G}(\mathbf{z})$. 

So now suppose that $(\mathbf{z},j(\mathbf{z}))\in V$ and $\dim_{G}(\mathbf{z}|C\cup\mathbf{t}) = \dim_{G}(\mathbf{z})$. Suppose that for some $1\leq i<k\leq n$ we have that there exists $g\in G$ such that $z_{i} = gz_{k}$. Then the tuple $(z_{1},\ldots,z_{n},a)$ is a solution of the following Khovanskii system (recall Remark \ref{rem:kclalgcl}):
\begin{equation*}
    \begin{array}{ccc}
         X_{1} &=& j(X_{n}) + X_{n+1}  \\
         X_{2} &=& j(X_{1})\\
         X_{3} &=& j(X_{2}) \\
         &\vdots& \\
         X_{n} &=& j(X_{n-1})\\
         X_{i} &=& gX_{k}
    \end{array}.
\end{equation*}
However, as this system is defined over $\mathbb{Q}$, then by Proposition \ref{prop:kleqj} that would mean that $a\in C$, which is a contradiction. Therefore we must have that $\dim_{G}(\mathbf{z}|C\cup\mathbf{t}) = n$.
\end{proof}

\subsection{Modular Schanuel conjecture}
\label{subsec:msc}
The classical statement of Schanuel's conjecture for the complex exponential function is the following:

\begin{conjecture}[Schanuel's conjecture, {{\cite[p. 30--31]{Lang-tr}}}]
If $x_{1},\ldots,x_{n}\in\mathbb{C}$ are $\mathbb{Q}$-linearly independent, then 
\begin{equation*}
    \mathrm{tr.deg.}_{\mathbb{Q}}\left(x_{1},\ldots,x_{n},\exp(x_{1}),\ldots,\exp(x_{n})\right)\geq n.
\end{equation*}
\end{conjecture}

As shown in \cite[1.3 Corollaire]{bertolin}, Schanuel's conjecture follows from the the \emph{generalised period conjecture} of Grothendieck-Andr\'e (see \cite[\textsection 23.4.4]{andre} and \cite{bertolin} for statements). Using the results of \cite{bertolin}, one can also specialise the generalised period conjecture to get a statement about the $j$-function. This can be done in different ways, which we now illustrate.

First, we recall that \cite{bertolin} proves that another consequence of the generalised period conjecture is the following Conjecture \ref{conj:mc}. The notation of the statement requires some explanation (see \cite{diaz} for terminology around elliptic curves). 

Let $E_{1},\ldots,E_{n}$ be elliptic curves, pairwise non-isogenous. For $E_{\ell}$, let $\omega_{1\ell}$ and $\omega_{2\ell}$ be its periods, let $\eta_{1\ell}$ and $\eta_{2\ell}$ be its quasi-periods, let $\tau_{\ell} = \frac{\omega_{1\ell}}{\omega_{2\ell}}\in\mathbb{H}^{+}$ (so that $j\left(\tau_{\ell}\right)$ is the $j$-invariant of $E_{\ell}$), let $q_{\ell} = \exp(2\pi i\tau_{\ell})$, and let $k_{\ell} = \mathrm{End}\left(E_{\ell}\right)\otimes_{\mathbb{Z}}\mathbb{Q}$. 
\begin{conjecture}[{{\cite[1.4 Remarque: Conjecture modulaire]{bertolin}}}]
\label{conj:mc}
\begin{equation*}
\mathrm{tr.deg.}_{\mathbb{Q}}\mathbb{Q}\left(2\pi i,\left\{ q_{\ell}, j\left(\tau_{\ell}\right), \omega_{1\ell}, \omega_{2\ell}, \eta_{1\ell}, \eta_{2\ell}\right\}_{\ell=1}^{n}\right) \geq \mathrm{rank}\left<q_{\ell}\right>_{\ell=1}^{n} + 4\sum_{\ell=1}^{n}\left(\dim_{\mathbb{Q}}k_{\ell}\right)^{-1} - n + 1,
\end{equation*}
where $\dim_{\mathbb{Q}}k_{\ell}$ is the dimension of $k_{\ell}$ as a $\mathbb{Q}$-vector space (i.e. it is the degree of the extension $k_{\ell}/\mathbb{Q}$, and so it is either 2 or 1, depending on whether $E_{\ell}$ is  CM or not, respectively), and $\mathrm{rank}\left<q_{\ell}\right>_{\ell=1}^{n}$ denotes the rank of the abelian multiplicative group generated by $q_{1},\ldots,q_{n}$.
\end{conjecture}

One can use the results conveniently compiled in \cite{diaz} to show that, for every $\tau\in\mathbb{H}^{+}$ satisfying $j'(\tau)\neq 0$, the following two fields have the same algebraic closure in $\mathbb{C}$ (we recall that $\eta_{1}\omega_{2}-\eta_{2}\omega_{1}=2\pi i$):
\begin{equation*}
    \mathbb{Q}\left(2\pi i, \exp(2\pi i \tau), j(\tau), \omega_{1}, \omega_{2}, \eta_{1}, \eta_{2}\right)\,\mbox{ and }\,
    \mathbb{Q}\left(2\pi i, \exp(2 \pi i \tau), \tau, j(\tau), j'(\tau), j''(\tau)\right).
\end{equation*}
Therefore, we can deduce the following statement from Conjecture \ref{conj:mc}.

\begin{conjecture}[Modular Schanuel conjecture with derivatives and special points]
Let $z_{1},\ldots,z_{n}\in \mathbb{H}^{+}$ be $G\mathrm{cl}$-independent and suppose that $j'(z_{i})\neq 0$, for $i=1,\ldots,n$. Let $m_{1} = \dim_{G}\left(\mathbf{z}|\Sigma\right)$, and let $m_{2} = n-m_{1}$. Then:
\begin{equation*}
    \mathrm{tr.deg.}_{\mathbb{Q}}\mathbb{Q}\left(\mathbf{z},j(\mathbf{z}),j'(\mathbf{z}),j''(\mathbf{z})\right)\geq 3m_{1} + m_{2}.
\end{equation*}
\end{conjecture}

This conjecture, along with the equivalence between (M1) and (M2), imply the following weaker statement, where we remove the presence of special points (see also \cite[Conjecture 8.3]{Pila-functional-transcendence}). 

\begin{conjecture}[Modular Schanuel conjecture with derivatives: MSCD]
Let $z_{1},\ldots,z_{n}\in \mathbb{H}^{+}$. Then:
\begin{equation*}
    \mathrm{tr.deg.}_{\mathbb{Q}}\mathbb{Q}\left(\mathbf{z},j(\mathbf{z}),j'(\mathbf{z}),j''(\mathbf{z})\right)\geq 3\dim_{G}(\mathbf{z}|\Sigma).
\end{equation*}
\end{conjecture}

\begin{remark}
\label{rem:mscd}
Saying that the $j$-field $(\C,\h)$ satisfies MSCD is equivalent to saying that $(S,\Sigma)\lhd (\C,\h)$, where $(S,\Sigma)$ is the $j$-subfield of $\C$ generated by the set $\Sigma$ of special points. 
\end{remark}

We can further remove the presence of derivatives and get a simpler statement (see also \cite[Conjecture 8.4]{Pila-functional-transcendence}). 

\begin{conjecture}[Modular Schanuel conjecture: MSC]
Let $z_{1},\ldots,z_{n}\in \mathbb{H}^{+}$. Then:
\begin{equation*}
    \mathrm{tr.deg.}_{\mathbb{Q}}\mathbb{Q}\left(\mathbf{z},j\left(\mathbf{z}\right)\right)\geq \dim_{G}(\mathbf{z}|\Sigma).
\end{equation*}
\end{conjecture}

\subsection{Conditional results}
\label{subsec:Jfinfields}
In this section we will use Theorem \ref{thm:extender3} to show that MSCD implies some analogous results to the ones in \textsection \ref{sec:genJfinfields}, but this time for subfields of $C=j\mathrm{cl}(\emptyset)$. We begin with an example showing how to get a version of Theorem \ref{thm:mainJ}.

\begin{example}
\label{ex:Jfinfield}
Let $t\in C$ be transcendental over $\mathbb{Q}$. By Theorem \ref{thm:extender3} we know that $t\in \mathrm{kcl}(\emptyset)$, which means that there exist $t=t_{1},\ldots,t_{n}\in C$ and $f_{1},\ldots,f_{n}\in\mathbb{Z}[X_{1},\ldots,X_{n}]^{J}$ satisfying a Khovanskii system. Furthermore, as explained in Remark \ref{rem:kclalgcl}, it is actually more convenient to choose our $j$-polynomials $f_{1},\ldots,f_{n}$ to be elements of the ring:
$$\overline{\mathbb{Q}}\left[\left\{g\mathbf{X},J(g\mathbf{X})\right\}_{g\in G}\right].$$
As the $j$-field $(\C,\h)$ is full, we can assume that $t_{1},\ldots,t_{n}\in\h$ (see Remark \ref{rem:kclsolinD}). Choose the system with $n$ minimal. As $\mathbf{t}$ is a non-singular solution of the Khovanskii system, MSCD implies that:
$$3\dim_{G}(\mathbf{t}|\Sigma)\leq \mathrm{tr.deg.}_{\mathbb{Q}}\mathbb{Q}(\mathbf{t},J(\mathbf{t})) \leq 3n.$$
By the minimality of $n$, we also get that the coordinates of $\mathbf{t}$ are $G\mathrm{cl}$-independent. Furthermore, any special point is in $\overline{\mathbb{Q}}\cap \h$, so, again by minimality of $n$, we may assume that no coordinate of $\mathbf{t}$ is special. Thus, MSCD implies:
$$\mathrm{tr.deg.}_{\mathbb{Q}}\mathbb{Q}(\mathbf{t},J(\mathbf{t})) = 3n.$$
\end{example}

\begin{theorem}
\label{thm:mainJ2}
Let $F$ be a subfield of $C$ such that $\mathrm{tr.deg.}_{\mathbb{Q}}F$ is finite. Then MSCD implies that there exist  $t_{1},\ldots,t_{m}\in\mathbb{H}$ such that:
\begin{enumerate}
    \item[(B1):] $F\subseteq\overline{\mathbb{Q}\left(\mathbf{t},J\left(\mathbf{t}\right)\right)}$,
    \item[(B2):] $\mathrm{tr.deg.}_{\mathbb{Q}}\mathbb{Q}\left(\mathbf{t},J\left(\mathbf{t}\right)\right) = 3\dim_{G}\left(\mathbf{t}\right|\Sigma)$.
\end{enumerate}
\end{theorem}
\begin{proof}
 If $F\subset\overline{\mathbb{Q}}$, then the result is trivial. So suppose $F$ has positive transcendence degree over $\mathbb{Q}$ and let $T$ be a transcendence basis for $\overline{F}$ such that $T\subset\h$. The idea now is to do something very similar to the construction in Example \ref{ex:Jfinfield}. By Theorem \ref{thm:extender3} we know that for every $t\in T$, $t\in \mathrm{kcl}(\emptyset)$, and so by extending Example \ref{ex:Jfinfield}, there exist $t_{1},\ldots,t_{n}\in D_{C}$ such that $T\subset\left\{t_{1},\ldots,t_{n}\right\}$ and $t_{1},\ldots,t_{n}$ are the solution of a Khovanskii system of $j$-polynomials. By choosing $n$ minimal, we can assume that the set $\left\{t_{1},\ldots,t_{n}\right\}$ is $G\mathrm{cl}$-independent (the elements of $T$ are clearly $G\mathrm{cl}$-independent) and that no element is a special point. Therefore 
 $$\mathrm{tr.deg.}_{\mathbb{Q}}\mathbb{Q}(\mathbf{t},J(\mathbf{t})) = 3n.$$
\end{proof}

Recall that 
$$\forall z\in\h\left(z\in\Sigma\iff\mathrm{tr.deg.}_{\mathbb{Q}}\mathbb{Q}(z,j(z)) =0\right).$$
Combined with the equivalence between (M1) and (M2), this implies that without loss of generality we can assume that in the statement of Theorem \ref{thm:mainJ2} the tuple $\mathbf{t}$ is $G\mathrm{cl}$-independent and none of its coordinates is special. In other words, we can assume that:
\begin{enumerate}
    \item[(B3):] $\dim_{G}\left(\mathbf{t}\right|\Sigma) = m$.  
\end{enumerate}

\begin{lemma}
\label{lem::mscdJf}
Suppose $t_{1},\ldots,t_{m}\in\h$ satisfy (B2). Then MSCD implies that for any $z_{1},\ldots,z_{n}\in \mathbb{H}$ we have:
$$\mathrm{tr.deg.}_{\mathbb{Q}(\mathbf{t},J(\mathbf{t}))}\mathbb{Q}(\mathbf{z},\mathbf{t},J(\mathbf{z}),J(\mathbf{t}))\geq3\dim_{G}(\mathbf{z}|\Sigma\cup\mathbf{t}).$$
\end{lemma}
\begin{proof}
The inequality is obtained by first using MSCD to get:
$$\mathrm{tr.deg.}_{\mathbb{Q}}\mathbb{Q}(\mathbf{z},\mathbf{t},J(\mathbf{z}),J(\mathbf{t}))\geq3\dim_{G}(\mathbf{z}\cup\mathbf{t}|\Sigma)$$
and now using the addition formula and (B2). 
\end{proof}

Just like with MSC, we can remove the presence of derivatives from Lemma \ref{lem::mscdJf} to get:

\begin{corollary}
\label{cor:Jfin}
Suppose $t_{1},\ldots,t_{m}\in\h$ satisfy (B2). Then MSCD implies that for any $z_{1},\ldots,z_{n}\in \mathbb{H}$ we have:
$$\mathrm{tr.deg.}_{\mathbb{Q}(\mathbf{t},J(\mathbf{t}))}\mathbb{Q}(\mathbf{z},\mathbf{t},j(\mathbf{z}),J(\mathbf{t}))\geq\dim_{G}(\mathbf{z}|\Sigma\cup\mathbf{t}).$$
\end{corollary}

\begin{remark}
\label{rem:notinC}
Assuming Conjecture \ref{conj:mc}, we can show that $\pi$ is not in $C$. For suppose that $\pi\in C$. Then by Theorem \ref{thm:mainJ2} (which can be used as Conjecture \ref{conj:mc} implies MSCD) we get that there exist $t_{1},\ldots,t_{n}\in\h$ such that $\pi$ is in the algebraic closure of $\mathbb{Q}(\mathbf{t},J(\mathbf{t}))$ and conditions (B2) and (B3) are satisfied. In the inequality of Conjecture \ref{conj:mc} we can eliminate the presence of the $q_{\ell}$ on the left-hand side if we also eliminate the rank of $\left<q_{\ell}\right>$ on the right-hand side, which gives us that
$$\mathrm{tr.deg.}_{\mathbb{Q}}\mathbb{Q}(2\pi i,\mathbf{t},J(\mathbf{t}))\geq 3n + 1,$$
(recall that condition (B3) ensures that no $t_{i}$ is special and that the elliptic curves associated to $t_{1},\ldots,t_{n}$ are non CM and pairwise non-isogenous). This implies then that
$$\mathrm{tr.deg.}_{\mathbb{Q}(\mathbf{t},J(\mathbf{t}))}\mathbb{Q}(2\pi i,\mathbf{t},J(\mathbf{t}))\geq 1,$$
thus showing that $\pi$ is transcendental over $\mathbb{Q}(\mathbf{t},J(\mathbf{t}))$, which is a contradiction. 
\end{remark}

\section*{Acknowledgements}
We are grateful to the referee for useful comments that helped us improve the presentation of the paper.

\bibliographystyle {alpha}
\bibliography {ref}

\end{document}